\documentclass[10pt]{amsart}

\usepackage{latexsym}
\usepackage{comment}
\usepackage{cancel} 

\usepackage{amsmath}
\usepackage{amsfonts}
\usepackage{amssymb}
 
\usepackage{graphicx,color}
\usepackage{epstopdf}
\usepackage{caption}
\usepackage{diagbox}
\usepackage{subfig}
\usepackage{bm}
\usepackage{algorithm}
\usepackage{algorithmic}
\usepackage{hyperref}
\usepackage{cleveref}
\usepackage{booktabs}
\usepackage{array}
\usepackage{enumerate}
\usepackage{tikz}
\usepackage{scalefnt}
\usepackage{hyperref}
\usepackage[margin=1in]{geometry}
\graphicspath{ {./figures/} }
\newtheorem{theorem}{Theorem}[section]

\theoremstyle{definition}

\newtheorem{example}[theorem]{Example}

\theoremstyle{remark}
\newtheorem{remark}[theorem]{Remark}

\numberwithin{equation}{section}

\newcommand{\bx}{{\bf{x}}}

\def\bz{\bm{\zeta}}

\def\O{\Omega}

\def\H2{{H^2(\O)}}

\def\eps{\varepsilon}

\DeclareMathOperator*{\argmin}{argmin}

\allowdisplaybreaks[1]

\begin{document}

\title[PINNs-OCP]{Two-scale neural networks for optimal control of linear convection-dominated equations}

\author{Sijing Liu}
\address{Sijing Liu, Department of Mathematical Sciences, Worcester Polytechnic Institute, 100 Institute Road, Worcester, MA 01609,
Current address: Department of Mathematics\\
SUNY Polytechnic Institute\\
Utica, NY\\
USA}
\email{sliu13@wpi.edu, lius7@sunypoly.edu}
\author{Marcus Sarkis}
\address{Marcus Sarkis, Department of Mathematical Sciences, Worcester Polytechnic Institute, 100 Institute Road, Worcester, MA 01609}
\email{msarkis@wpi.edu}

\author{Yi Zhang}
\address{Yi Zhang, Department of Mathematics and Statistics, The University of North Carolina at Greensboro, 1400 Spring Garden Street, Greensboro, NC 27412}
\email{y\_zhang7@uncg.edu}

\author{Zhongqiang Zhang}
\address{Zhongqiang Zhang, Department of Mathematical Sciences, Worcester Polytechnic Institute, 100 Institute Road, Worcester, MA 01609}
\email{zzhang7@wpi.edu}
\date{\today} 

\begin{abstract}
We propose a two-scale neural network method for optimal control problems governed by convection‑dominated convection-diffusion-reaction equations. Building on two‑scale architectures developed for singularly perturbed forward problems, we augment the spatial input with suitably rescaled features that become increasingly important as the diffusion 
coefficient 
becomes small. 
The approach employs separate neural networks for the state and adjoint state variables of the optimality system, reflecting the fact that these quantities develop sharp layers in different parts of the domain due to opposite convection fields. 
By choosing different center points for the two networks, the architecture naturally aligns with the layer location of each variable.
We present two formulations of the method, one based on the first‑order optimality conditions and another using penalization of the PDE constraint, and combine them with a successive training strategy that gradually decreases the diffusion coefficient toward its target value. 
Numerical experiments on benchmark problems illustrate the effectiveness and behavior of the proposed approach.


\end{abstract} 

\maketitle


\section{Introduction}

Optimal control problems (OCPs) constrained by partial differential equations (PDEs) provide a mathematical framework for selecting inputs that steer a physical system toward a desired state while accounting for control cost. They arise in aerodynamics, shape optimization, biology, economics, and many other areas, and have been studied extensively in numerical analysis and optimization \cite{hinze2008optimization,Lions,Tro}. In a typical tracking-type formulation, the state satisfies a PDE constraint, and the objective balances the mismatch between the state and a target profile against a regularization term for the control.

We consider OCPs governed by convection-dominated convection-diffusion-reaction equations. When the diffusion parameter $\eps>0$ is small, the state equation is singularly perturbed: solutions may be smooth in much of the domain but develop thin boundary or interior layers with  width related to $\eps$. The optimal control setting introduces an additional structural difficulty. The first-order optimality system couples the state equation to an adjoint equation whose convection direction is reversed. Consequently, the state and adjoint can develop layers in different parts of the domain, often near opposite inflow/outflow boundaries. A numerical method must therefore approximate coupled variables with different layer locations while preserving the optimality relation among state, adjoint, and control.

The focus of this work is to explore whether explicitly enforcing the first-order optimality system improves physics-informed neural network (PINN) approximations for such singularly perturbed OCPs. 
In classical PDE-constrained optimization, adjoint equations encode first-order sensitivity information and provide the standard route to optimality conditions. In a PINN setting, however, one can either train networks against the residuals of the state-adjoint optimality system or avoid the adjoint equation by penalizing the PDE constraint directly in the objective. These choices lead to different residual structures, optimization landscapes, and layer-resolution demands. The convection-dominated regime provides a natural test case because the adjoint equation is not merely an additional equation; it carries layer information transported in the opposite direction from the state.

\subsection{Existing approaches and limitations}

Classical discretizations for convection-dominated equations are shaped by a balance among stability, layer resolution, and interior accuracy. Standard Galerkin finite element methods on quasi-uniform meshes may exhibit spurious oscillations near sharp layers, motivating stabilized schemes such as streamline-upwind/Petrov--Galerkin methods \cite{Hughes1979MULTIDIMENSIONALUS,BROOKS1982199}, edge-average finite element methods \cite{adler2023stable,xu1999monotone}, local projection stabilization \cite{knobloch2009local}, and discontinuous Galerkin methods \cite{doi:10.1137/S0036142997316712,doi:10.1137/S0036142900374111}. 
For OCPs, the discretization must also be consistent with the state-adjoint-control structure, leading to the well-known distinction between optimize-then-discretize and discretize-then-optimize approaches. Weak imposition of boundary conditions, as in discontinuous Galerkin methods \cite{leykekhman2012local,liu2024multigrid} and Nitsche's method \cite{nitsche1971variationsprinzip}, can improve interior accuracy, but physical boundary layers may be ignored when $\eps$ is far less than the mesh size  \cite{leykekhman2012local}. HDG methods for convection-dominated OCPs have also been developed \cite{chen2018hdg,liu2025balancing,liu2025convergence}. These results highlight a persistent tradeoff: choices that stabilize the computation or improve interior behavior do not necessarily resolve the singular layers themselves.

Deep learning methods provide a mesh-free alternative, but they do not remove the layer-resolution difficulty. PINN-type methods for layer problems can be viewed broadly as either layer-informed or layer-agnostic. Layer-informed approaches use explicit knowledge of the layer location or type \cite{gao2025more,wang2024general,wang2024aspinn}. Layer-agnostic approaches instead seek mechanisms for representing sharp gradients without prescribing where they occur, including curriculum learning \cite{cao2024multistep,cao2023physics,qiao2025two}, adaptive loss weighting \cite{wang2024less}, and architectures designed for steep transitions \cite{gao2025more}. The present work does not prescribe an analytic layer profile or exact layer width; instead, it incorporates available inflow/outflow information through layer-informed center placement.

For deep learning approaches to OCPs, an equally important distinction is whether adjoint equations are enforced explicitly. Non-adjoint formulations impose the PDE constraint through penalties or related constrained-learning mechanisms, including adversarial and adaptive-weight strategies \cite{cao2025adversarial} and value-function/control approximation methods \cite{dupret2026deep}. Adjoint-based formulations incorporate first-order optimality information directly, with applications to error estimation \cite{dai2025solving}, PINN stability \cite{Nzoyem2023comparison}, shape optimization \cite{wang2024aonn,yin2024aonn}, and staged training \cite{barry2025physics}; see also the review \cite{wang2026optimal}. What remains unclear is how this formulation choice affects PINN training for convection-dominated OCPs, where the state and adjoint equations transport layers in opposite directions. This gap motivates the comparison in this paper.

\subsection{Current work: contributions and significance}

The paper compares two PINN formulations for convection-dominated OCPs within a common two-scale continuation framework. The optimality-condition formulation approximates the state and adjoint with separate networks and trains them against the residuals of the first-order optimality system. The penalization formulation enforces the PDE constraint through a penalty term and does not introduce the adjoint equation explicitly. Placing these formulations in the same computational framework isolates the effect of adjoint enforcement in a singularly perturbed training problem.

The common framework is based on the two-scale construction of \cite{qiao2025two}. The network input is augmented with rescaled features depending on the diffusion coefficient $\eps$, and training proceeds by decreasing $\eps$ toward the target value. This continuation strategy is compatible with relatively lightweight collocation point distributions \cite{munzer2022curriculum}. Separate two-scale networks are used for variables whose layers may occur at different locations, the state and adjoint in the optimality-condition formulation, and the state and control in the penalization formulation. Each variable is therefore represented using its own layer-informed center placement and rescaled coordinates, without prescribing the precise layer location or width.

The main contribution lies in the formulation-level study of PINNs for singularly perturbed optimality systems. The comparison addresses whether enforcing the adjoint equation supplies useful training structure beyond what is obtained from penalizing the PDE constraint. This question is specific to OCPs and becomes especially relevant in the convection-dominated regime, where the adjoint carries layer information distinct from that of the state.

The contributions are summarized as follows:
\begin{itemize}
    \item[(i)] A PINN framework is developed for convection-dominated OCPs in two formulations: an optimality-condition formulation that enforces the state equation, adjoint equation, and optimality, and a penalization formulation that avoids the explicit adjoint equation.

    \item[(ii)] A controlled comparison is carried out between adjoint-based and non-adjoint PINN formulations for singularly perturbed OCPs whose state and adjoint equations have opposite convection directions and distinct layer structures.

    \item[(iii)] The two-scale architecture and continuation strategy in \cite{qiao2025two} are inherited and adapted here from prior two-scale PINN work to the setting of coupled optimality systems. While the underlying two-scale construction is based on these existing approaches, our method applies it to optimal control problems, where the coupled state, adjoint, and control variables exhibit a different underlying structure. This setting is distinct from the application of two-scale neural networks to systems of ordinary differential equations \cite{qiao2026two}.

    \item[(iv)] A successive training strategy in the diffusion parameter is adapted to the OCP setting. Starting from a moderate value $\eps_0$, the diffusion coefficient is decreased toward the target value, with each stage warm-started from the previous one.
\end{itemize}

The numerical results indicate that, within this common two-scale setting, explicit enforcement of the first-order optimality system yields smaller residuals and improved accuracy compared with penalization-based training. The results support the view that, for convection-dominated OCPs, the adjoint equation provides useful structure for PINN training rather than serving only as an analytical optimality device.

\subsection{Outline of the paper}

The remainder of the paper is organized as follows. Section \ref{sec:ocp} presents the OCP under consideration, derives the first-order optimality system, and introduces the penalization reformulation. Section \ref{sec:PINN} reviews physics-informed neural networks and the two-scale neural network architecture. Section \ref{sec:double-PINNs} presents the proposed two-scale neural network method for OCPs in both the optimality-condition and penalization formulations, together with the successive training algorithm for small $\eps$. Section \ref{sec:numer} reports numerical experiments on benchmark problems. Section \ref{sec:summary} concludes with a summary and discussion.

\section{An optimal control problem}\label{sec:ocp}
This section considers a linear-quadratic optimal control problem constrained by a convection-dominated convection-diffusion-reaction equation. This model problem admits two numerical difficulties that are central in singularly perturbed optimal control: the state and adjoint may develop layers in different parts of the domain, and the choice of PINN loss determines how these layers enter the training problem:
\begin{equation}\label{optcon}
(\bar{y},\bar{u})=\argmin_{(y,u)} J(y,u), 
\end{equation} 
subject to the state equation 
\begin{equation}\label{eq:cdr-intro}
  -\eps\Delta y + \bz\cdot\nabla y + c y = f + u \quad \text{in } \Omega,
  \qquad y = 0 \quad \text{on } \partial \Omega, 
\end{equation} 
with the cost functional 
\begin{align}
     J(y,u)&=\frac{1}{2}\|y-y_d\|^2_{L^2(\Omega)}+\frac{\beta}{2}\|u\|^2_{L^2(\Omega)}. 
\end{align}
Here $\eps>0$ is the diffusion parameter, $y_d\in L^2(\Omega)$ is the desired state, $f\in L^2(\Omega)$ is a source term, $\bm{\zeta}\in [W^{1,\infty}(\Omega)]^n$ is the convective field, $c \in L_{\infty}(\Omega)$ is nonnegative, $\beta > 0$ is a regularization constant, and $\Omega \subset \mathbb{R}^2$ is a bounded open polygonal domain with Lipschitz boundary $\partial \Omega$. 
The state and control are sought in $H^1_0(\Omega)$ and $L^2(\Omega)$, respectively. Throughout, we assume 
\begin{equation}\label{eq:advassump}
    c-\frac12\nabla\cdot\bz\ge c_0>0
\end{equation}
for some constant $c_0$, which ensures \eqref{eq:cdr-intro} is well-posed \cite{ayuso2009discontinuous,liu2024multigrid}. 

When $\eps \ll |\bz|$, solutions may contain sharp outflow boundary layers and, depending on the data and geometry, interior layers aligned with characteristic curves of $\bz$. The layer widths are problem dependent and vanish as $\eps \to 0$, so standard PINN training can under-resolve the regions that dominate the residual.

\subsection{First-order optimality condition} \label{subsec:OC}
The minimizer of \eqref{optcon} is characterized by the first-order optimality system \cite{Lions, Tro}. Introducing an adjoint state $p \in H^1_0(\Omega)$ gives
\begin{subequations}\label{eq:sp}
\begin{alignat}{3}
-\eps\Delta p-\bz\cdot\nabla p+(c-\nabla\cdot\bm{\zeta})p &= y - y_d, \quad &&\text{in } \Omega, \qquad p = 0 \quad \text{on } \partial \Omega, \label{eq:adjoint} \\
p+\beta u&=0, \quad &&\text{in } \Omega, \label{eq:spu}\\  
 -\eps\Delta y + \bz\cdot\nabla y + c y &= f + u \quad &&\text{in } \Omega,
  \qquad y = 0 \quad \text{on } \partial \Omega. \label{eq:state-3}
\end{alignat}
\end{subequations}
The adjoint equation \eqref{eq:adjoint} is again a convection-diffusion-reaction equation, but its convection field is $- \bm{\zeta}$ rather than $+ \bm{\zeta}$. 
Consequently, the adjoint develops layers in locations dual to those of the state. For PINN discretizations this dual-layer structure is not only a resolution issue; it also affects whether the loss explicitly enforces the adjoint equation or recovers the control through penalization. This observation motivates the formulation comparison and the double two-scale architecture introduced below.

Eliminating $u$ from \eqref{eq:spu} gives the coupled state-adjoint system  
\begin{subequations}\label{eq:osp}
\begin{alignat}{2}
-\eps\Delta p-\bz\cdot\nabla p+(c-\nabla\cdot\bm{\zeta})p &= y - y_d, \quad \text{in } \Omega, \qquad p = 0 \quad \text{on } \partial \Omega, \label{eq:adjoint-2} \\
 -\eps\Delta y + \bz\cdot\nabla y + c y + \frac{1}{\beta} p &= f \quad \text{in } \Omega,
  \qquad y = 0 \quad \text{on } \partial \Omega. 
\end{alignat}
\end{subequations} 

Define the operator $\mathcal{L}: H^1_0(\Omega) \longrightarrow H^{-1}(\Omega)$ by 
\[
\mathcal{L} y = -\eps\Delta y + \bz\cdot\nabla y + c y.
\]
Its formal adjoint in $L^2(\Omega)$ with respect to the standard inner product is 
\[
\mathcal{L}^* p = -\eps\Delta p - \bz\cdot\nabla p + (c-\nabla\cdot\bm{\zeta}) p, 
\]
so \eqref{eq:osp} reads 
\begin{subequations}\label{eq:speps1}
\begin{alignat}{2}
\mathcal{L}^* p - y &= - y_d, \quad \text{in } \Omega, \qquad p = 0 \quad \text{on } \partial \Omega, \label{eq:adjoint-op} \\
 \mathcal{L} y + \frac{1}{\beta} p &= f \quad \text{in } \Omega,
  \qquad y = 0 \quad \text{on } \partial \Omega. 
\end{alignat}
\end{subequations}

\subsection{Penalization} \label{subsec:PM}
An alternative route is to approximate \eqref{optcon} by penalization (cf. \cite[Chapter 5, Section 3]{Lions} and \cite{bergounioux1992penalization}). This formulation is attractive for PINNs because it avoids writing the adjoint equation explicitly, but it shifts the numerical burden to penalty parameter selection. Let $\alpha=\{\alpha_1,\alpha_2\}$ with $\alpha_i>0$ for $i=1,2$, and define the penalized cost
\begin{equation}\label{eq:penalization}
    J_\alpha(y,u)=J(y,u)+\alpha_1\|\mathcal{L}y-u-f\|^2_{L^2(\Omega)}+\alpha_2\|y\|^2_{L^2(\partial\Omega)},
\end{equation}
and the penalized problem
\begin{equation}\label{optconp}
(\bar{y}_\alpha,\bar{u}_\alpha)=\argmin_{(y,u) \in H^2(\Omega)\times L^2(\Omega)} J_\alpha(y,u).
\end{equation}
In \eqref{eq:penalization}, the first penalty measures the strong-form PDE residual in $\Omega$, while the second measures violation of the homogeneous Dirichlet condition on $\partial \Omega$. Thus the constrained problem is replaced by an unconstrained minimization over $(y,u)$, a convenient form for neural network training but one whose accuracy depends on the penalty weights.

The following theorem is a variant of the results in \cite[Chapter 5, Section 3]{Lions} and \cite[Theorem 3.3]{bergounioux1992penalization}.
\begin{theorem}\label{thm:pen}
    As $\alpha=\{\alpha_1, \alpha_2\}\rightarrow\infty$, we have
    \begin{equation}
        \min_{(y,u) \in H^2(\Omega)\times L^2(\Omega)} J_\alpha(y,u)\rightarrow\min_{(y,u) \in H^1_0(\Omega)\times L^2(\Omega)} J(y,u).
    \end{equation}
Moreover,
    \begin{equation}
    \bar{y}_\alpha\rightarrow\bar{y}\quad\text{in}\ L^2(\Omega), \qquad  
    \bar{u}_\alpha\rightarrow\bar{u}\quad\text{in}\ L^2(\Omega).
    \end{equation}
\end{theorem}
\begin{remark}
     The penalization formulation avoids explicit enforcement of the adjoint equation and extends naturally to nonlinear state equations (cf. \cite{Lions}).
\end{remark}

\section{Physics-informed neural networks}\label{sec:PINN}
This section recalls the neural-network and PINN notation used later, then reviews the two-scale PINN architecture of \cite{qiao2025two}. The emphasis is on the features that are needed for singular perturbations: residual evaluation by automatic differentiation and input enrichment by a stretched coordinate.

\subsection{Neural Networks and PINNs} 

A feedforward neural network of $L$ layers with input $x \in \mathbb{R}^{n_0}$ is defined recursively as
\begin{subequations}\label{eq:nnstruc}
\begin{alignat}{2}
    a^0(\bm{x})&=\bm{x}\\
    a^l(\bm{x})&=\sigma(W^la^{l-1}(\bm{x})+b^l),\quad 1\le l<L, \\
    a_{N\!N}(\bm{x})&=W^L a^{L-1}(\bm{x})+b^L 
\end{alignat}
\end{subequations}
where $W^l\in\mathbb{R}^{n_l\times n_{l-1}}$ and $b^l\in\mathbb{R}^{n_l}$ are the weights and biases of layer $l$, and $\sigma$ is a nonlinear activation function applied componentwise (e.g, tanh or rectified linear unit (ReLU)). 
We write $\mathcal{N}_{\theta}(\bm{x}) := a_{NN}(x)$ with parameters $\theta = \{ W^\ell, b^\ell \}_{\ell=1}^{L}$.

Given a boundary value problem 
\begin{equation}\label{eq:dfq}
        \mathcal{D}w=g(x)\quad \mbox{in}\ \Omega, \qquad \mathcal{B}w=h(x)\quad \mbox{on}\ \partial\Omega,
\end{equation}
with differential operators $\mathcal{D}$ and $\mathcal{B}$, the PINN approximation $w_{\theta}(x) := \mathcal{N}_{\theta}(x)$ is obtained by minimizing the loss 
\begin{equation}\label{eq:jd}
    J_d(\theta)=\frac{1}{N_c}\sum_{k=1}^{N_c}|g(\bm{x}_r^k)-\mathcal{D}w_\theta (\bm{x}_r^k)|^2+\frac{\alpha}{N_b}\sum_{k=1}^{N_b}|h(\bm{x}_b^k)-\mathcal{B}w_\theta (\bm{x}_b^k)|^2,
\end{equation}
where $\{\bm{x}_r^k\}_{k=1}^{N_c} \subset \Omega$ are interior collocation points, $\{\bm{x}_b^k\}_{k=1}^{N_b} \subset \partial \Omega$ are boundary collocation points, and $\alpha \ge 0$ is the weight balancing the PDE residual against the boundary loss.  
The differential operators are evaluated by automatic differentiation. We use the Adam optimizer \cite{kingma2017adammethodstochasticoptimization} to minimize $J_d(\theta)$.

 
\subsection{Two-scale PINNs}
For singularly perturbed PDEs, the approximation space should represent both the outer solution and rapid variation inside thin layers. The two-scale neural network of \cite{qiao2025two} addresses this by augmenting the physical input $x$ with two additional features that encode the small parameter $\eps$. 

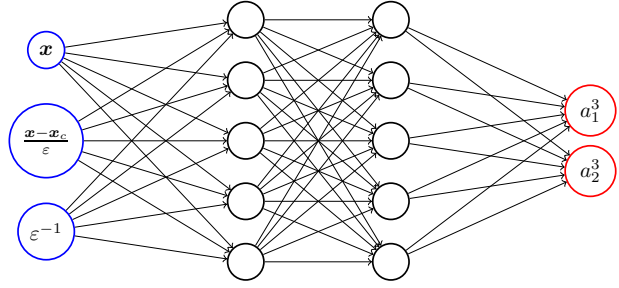
\begin{figure}[h!] 
    \centering
    \caption{Two-scale neural network}
\label{figure:tnnocp}
    \scalebox{0.8}{
\begin{tikzpicture}[
    input/.style={circle, draw=blue, thick, minimum size=6mm},
    hidden/.style={circle, draw, thick, minimum size=6mm},
    output/.style={circle, draw=red, thick, minimum size=6mm},
    dots/.style={minimum size=0mm, draw=none, fill=none},
    layerlabel/.style={rectangle, draw=none, text centered}
]

\node[input] (I1) at (-1.5,1.5) {$\bm{x}$};
\node[input] (I2) at (-1.5,0) {$\frac{\bm{x}-\bm{x}_c}{\eps}$};
\node[input] (I3) at (-1.5,-1.5) {$\eps^{-1}$};

\node[hidden] (H11) at (1.8,2) {};
\node[hidden] (H12) at (1.8,1) {};
\node[hidden] (H13)at (1.8,0) {};
\node[hidden] (H14) at (1.8,-1) {};
\node[hidden] (H15) at (1.8,-2) {};

\node[hidden] (H21) at (4.2,2){};
\node[hidden] (H22) at (4.2,1){};
\node[hidden] (H23) at (4.2,0){};
\node[hidden] (H24) at (4.2,-1){};
\node[hidden] (H25) at (4.2,-2){};

\node[output] at (7.5,0.5) (O1) {$a^3_1$};
\node[output] at (7.5,-0.5) (O2) {$a^3_2$};

\foreach \i in {1,2,3}
  \foreach \j in {11,12,13,14,15}
    \draw[->] (I\i) -- (H\j);

\foreach \i in {11,12,13,14, 15}
  \foreach \j in {21,22,23,24, 25}
    \draw[->] (H\i) -- (H\j);
\foreach \i in {21,22,23,24,25}
  \foreach \j in {1,2}
    \draw[->] (H\i) -- (O\j);

\end{tikzpicture}}
\end{figure}

Let $\Omega \subset \mathbb{R}^d$ be the spatial domain and let $x_c \in \mathbb{R}^d$ be a chosen center (typically the geometric center of $\Omega$ or a point near the layer when the location of large gradients is known a priori). 
Let $\gamma < 0$ be a fixed scaling exponent. The two-scale network takes the augmented vector 
\[
(\bm{x}, \eps^\gamma(\bm{x}-\bm{x}_c), \eps^{\gamma}) \in \mathbb{R}^{2d+1} 
\]
as input, 
where $\eps^\gamma(\bm{x}-\bm{x}_c) \in \mathbb{R}^d$ is the rescaled coordinate around $x_c$, and $\eps^\gamma$ supplies the magnitude of the small parameter. The forward network is defined as in \eqref{eq:nnstruc}, with this enriched input:
\begin{subequations}\label{eq:tnnstruc}
\begin{alignat}{2}
    a^0(\bm{x},\eps^\gamma(\bm{x}-\bm{x}_c),\eps^{\gamma})&=[\bm{x}, \eps^\gamma(\bm{x}-\bm{x}_c), \eps^{\gamma}]\\
    a^l(\bm{x}{,\eps^\gamma(\bm{x}-\bm{x}_c),\eps^{\gamma}})
    &=\sigma(W^la^{l-1}(\bm{x}
    {,\eps^\gamma(\bm{x}-\bm{x}_c),\eps^{\gamma}})+b^l),\ 1 \leq l<L,  \\
    a_{N\!N}(\bm{x},\eps^\gamma(\bm{x}-\bm{x}_c),\eps^{\gamma})&=W^La^{L-1}(\bm{x},\eps^\gamma(\bm{x}-\bm{x}_c),\eps^{\gamma})+b^L.
\end{alignat}
\end{subequations}
We denote the resulting network by $\mathcal{N}_{\theta}(\bm{x},\eps^\gamma(\bm{x}-\bm{x}_c),\eps^{\gamma})$. The PINN loss is then defined as in \eqref{eq:jd} with $w_\theta(x) := \mathcal{N}_{\theta}(\bm{x},\eps^\gamma(\bm{x}-\bm{x}_c),\eps^{\gamma})$. 
Figure \ref{figure:tnnocp} shows the architecture with  $\gamma=-1$ and two output variables. 

The rescaled feature has two complementary effects. First, the chain rule gives 
\[
\frac{\partial}{\partial x_j} \mathcal{N}_{\theta} (\bm{x},\eps^\gamma(\bm{x}-\bm{x}_c),\eps^{\gamma}) = \frac{\partial \mathcal{N}_\theta}{\partial x_j} + \eps^\gamma \frac{\partial \mathcal{N}_\theta}{\partial (\eps^\gamma (x_j - x_{c,j}))},  
\]
where the first term is the derivative through the physical input and the second is the derivative through the rescaled input. With $\gamma < 0$ and $\eps \ll 1$, the factor $\eps^\gamma$ is large, so the network can produce derivatives of order $\eps^\gamma$ to match large derivatives that occur within boundary or interior layers of the singularly perturbed PDE. 
Second, by analogy with the matched method of asymptotic expansions, the solution near a layer admits an inner expansion in a stretched coordinate. For a layer at $\boldsymbol{x}_c$, the stretched coordinate has the form $(\boldsymbol{x} - \boldsymbol{x}_c) / \eps^\mu$ for some problem-dependent $\mu > 0$. In the convection-diffusion setting one typically has $\mu = \frac12$ or $\mu = 1$. The two-scale network exposes this stretched coordinate directly through its rescaled input $\eps^\gamma(\bm{x}-\bm{x}_c)$, with $\gamma = -\mu$. The appropriate combination of $\boldsymbol{x}$ and the stretched coordinate is learned during training, so no explicit asymptotic decomposition is required. 
The scalar feature $\eps^\gamma$ gives the network direct access to the value of $\eps$, independently of $\boldsymbol{x}$. This is useful when the solution has an $\eps$-dependent amplitude that does not vary spatially and during successive training, where $\eps$ changes between stages.

The center $\bm{x}_c$ determines where the rescaled coordinate $\eps^\gamma(\bm{x}-\bm{x}_c)$ is small. If the layer location is unknown, $\bm{x}_c$ may be chosen as the geometric center of $\Omega$, allowing the first hidden layer to learn an effective offset. If the layer location is known a priori, $\bm{x}_c$ can be placed near the layer so that the rescaled coordinate is small where the solution varies rapidly. Section \ref{sec:double-PINNs} uses this layer-informed choice separately for the state and adjoint/control networks. 
 
%
\section{Double Two-scale PINNs for Optimal Control Problems}\label{sec:double-PINNs}
This section presents the PINN methodology for \eqref{optcon}. The construction is organized around two choices that are consequential for singularly perturbed optimal control systems. The first is architectural: the state and the adjoint, or the state and the control, are represented by separate two-scale networks so that their stretched coordinates can be centered near different layer regions. The second is variational: training is based either on the first-order optimality system of Section \ref{subsec:OC} or on the penalized formulation of Section \ref{subsec:PM}. This organization makes the formulation choice explicit and separates it from the continuation strategy used for small~$\eps$.

 
\subsection{Double neural networks}
\label{subsec:double-nn}
 
The same architectural principle is used for both formulations. Although a single two-scale network with two outputs could represent both the state $y$ and the adjoint $p$ (or the control $u$), such a representation uses one center for variables whose layers generally occur in different locations. We therefore use separate two-scale networks:
\begin{align}
y_{\theta_y}(\bx) &= \mathcal{N}^y_{\theta_y}\!\bigl(\bx,\; \eps^\gamma(\bx - \bx_c^y),\; \eps^\gamma\bigr), \label{eq:y-net}\\
w_{\theta_w}(\bx) &= \mathcal{N}^w_{\theta_w}\!\bigl(\bx,\; \eps^\gamma(\bx - \bx_c^w),\; \eps^\gamma\bigr), \label{eq:w-net}
\end{align}
where $w = p$ in the optimality-condition formulation of Section \ref{subsec:OC} and $w = u$ in the penalization formulation of Section \ref{subsec:PM}. The two networks share the value of $\gamma$ but have different centers $\bx_c^y$ and $\bx_c^w$.
 
Separate networks with separate centers reflect the dual-layer structure identified in Section \ref{sec:ocp}. By \eqref{eq:spu}, the control satisfies $u = -p/\beta$ and therefore inherits the layer structure of $p$. The state equation \eqref{eq:state-3} has convection field $+\bz$, whereas the adjoint equation \eqref{eq:adjoint} has convection field $-\bz$. 
Following \cite{roos2008robust}, the layers of a convection-diffusion problem are determined by the convection field, and reversing the convection moves the layers to dual locations.  

A single two-scale network with one center cannot place the rescaled coordinate $\eps^\gamma(\bx - \bx_c)$ near both layer regions simultaneously. The double-network architecture addresses this mismatch directly: each network has its own center, chosen so that the rescaled coordinate is small near the layer of the corresponding variable. In the experiments we set
\begin{equation}
\bx_c^y \;=\; \bx_{\mathrm{outflow}}, \qquad
\bx_c^w \;=\; \bx_{\mathrm{inflow}}, 
\end{equation}
where $\bx_{\mathrm{outflow}}$ and $\bx_{\mathrm{inflow}}$ are points on the outflow and inflow boundaries of the state problem. 
 
The two networks have separate parameters but are trained through a coupled loss. With $L$ hidden layers of width $W$, the resulting architecture is $(2d+1, W, W, \ldots, W, 1)$ for each network. The parameter count is therefore twice that of a single shared network of the same depth and width, a cost that is offset by the ability to place the stretched coordinates at different layer locations. 

 
\subsection{Solving the first-order optimality system}
\label{subsec:solve-oc}
 
The first formulation trains directly on the reduced first-order optimality system \eqref{eq:osp}. This loss treats the state and adjoint equations symmetrically and recovers the control from $p+\beta u=0$ after training. 
 We define the loss function
\begin{equation}\label{eq:jodocp}
J^o_d(\theta_y, \theta_p) \;=\; R_{\mathrm{state}}(\theta_y, \theta_p) + R_{\mathrm{adj}}(\theta_y, \theta_p) + B_y(\theta_y) + B_p(\theta_p),
\end{equation}
where each term corresponds to the residual of a different equation in the optimality system, evaluated at collocation points. The interior residual terms are
\begin{align}
R_{\mathrm{state}}(\theta_y, \theta_p)
  &= \frac{1}{N_c} \sum_{k=1}^{N_c} \bigl| \mathcal{L} y_{\theta_y}(\bx_r^k) + \tfrac{1}{\beta}\, p_{\theta_p}(\bx_r^k) - f(\bx_r^k) \bigr|^2, \label{eq:Ra}\\
R_{\mathrm{adj}}(\theta_y, \theta_p)
  &= \frac{1}{N_c} \sum_{k=1}^{N_c} \bigl| \mathcal{L}^* p_{\theta_p}(\bx_r^k) - y_{\theta_y}(\bx_r^k) + y_d(\bx_r^k) \bigr|^2, \label{eq:Rb}
\end{align}
and the boundary terms are 
\begin{equation}\label{eq:Byp}
B_y(\theta_y) = \frac{\alpha_y}{N_b} \sum_{k=1}^{N_b} \bigl| y_{\theta_y}(\bx_b^k) \bigr|^2,
\qquad
B_p(\theta_p) = \frac{\alpha_p}{N_b} \sum_{k=1}^{N_b} \bigl| p_{\theta_p}(\bx_b^k) \bigr|^2,
\end{equation}
with weights $\alpha_y, \alpha_p > 0$ chosen to balance the boundary loss against the interior loss. 
Here $\{\bx_r^k\}_{k=1}^{N_c} \subset \Omega$ and $\{\bx_b^k\}_{k=1}^{N_b} \subset \partial\Omega$ are collocation points. 

The training problem is 
\begin{equation}
(\theta_y^\ast, \theta_p^\ast) \;=\; \operatorname*{argmin}_{(\theta_y, \theta_p)} \; J^o_d(\theta_y, \theta_p), \label{eq:OCtrain}
\end{equation}
and the optimal control is recovered as $u_{\theta_p^\ast} = -p_{\theta_p^\ast} / \beta$. The PDE residuals in \eqref{eq:Ra}--\eqref{eq:Rb} are evaluated pointwise via automatic differentiation. 


\subsection{Solving the penalized problem}
\label{subsec:solve-pm}
 
The second formulation trains on the penalized problem \eqref{optconp}. Since this problem involves the state $y$ and control $u$ directly, the two networks represent $y$ and {\color{cyan}$u$}. This loss is closer to an unconstrained optimization problem, but the PDE and boundary conditions are enforced through finite penalty weights. Note that no adjoint equation is involved in the penalized approach.
 
The loss function is a direct discretization of the penalized cost functional \eqref{eq:penalization}:
\begin{equation}\label{eq:jpdocp} 
J^p_d(\theta_y, \theta_u) 
  \;=\; T(\theta_y, \theta_u) + \alpha_1 \, R_{\mathrm{pen}}(\theta_y, \theta_u) + \alpha_2 \, B_y(\theta_y),
\end{equation}
where the tracking term, the PDE-residual penalty term, and the boundary penalty term are
\begin{align}
T(\theta_y, \theta_u)
  &= \frac{1}{2 N_c} \sum_{k=1}^{N_c} \Bigl[ \bigl| y_{\theta_y}(\bx_r^k) - y_d(\bx_r^k) \bigr|^2 + \beta \bigl| u_{\theta_u}(\bx_r^k) \bigr|^2 \Bigr], \label{eq:Jpd:a}\\
R_{\mathrm{pen}}(\theta_y, \theta_u)
  &= \frac{1}{N_c} \sum_{k=1}^{N_c} \bigl| \mathcal{L} y_{\theta_y}(\bx_r^k) - u_{\theta_u}(\bx_r^k) - f(\bx_r^k) \bigr|^2, \label{eq:Jpd:b}\\
B_y(\theta_y)
  &= \frac{1}{N_b} \sum_{k=1}^{N_b} \bigl| y_{\theta_y}(\bx_b^k) \bigr|^2. \label{eq:Jpd:c}
\end{align}
The training problem is
\begin{equation}
(\theta_y^\ast, \theta_u^\ast) \;=\; \operatorname*{argmin}_{(\theta_y, \theta_u)} \; J^p_d(\theta_y, \theta_u). \label{eq:PMtrain}
\end{equation}
Similar to Theorem \ref{thm:pen}, as $\alpha_1, \alpha_2 \to \infty$ the solution of \eqref{eq:PMtrain} may approach the solution of the original constrained problem \eqref{optcon} under appropriate assumptions on the continuous problem and the collocation points in the loss function, e.g. \cite{jin2023solving}. However, this continuous penalty argument does not by itself establish convergence of the finite-collocation PINN training procedure. The latter involves additional considerations, including sampling/quadrature consistency, neural approximation capacity, and nonconvex optimization. In practice, $\alpha_1$ and $\alpha_2$ are chosen large but finite, balancing the strength of the constraint enforcement against numerical conditioning of the loss. Thus, convergence of the trained collocation PINN should be viewed as an empirical numerical question unless additional assumptions are established.
 
The two formulations therefore test different ways of incorporating the optimal control structure into a PINN. The optimality formulation enforces the adjoint equation explicitly and recovers the control by post-processing, at the cost of a coupled state-adjoint residual. The penalized formulation avoids the adjoint equation and is structurally simpler, but its accuracy depends on penalty weights and conditioning. The numerical experiments in Section \ref{sec:numer} compare these tradeoffs.


\subsection{Training the double neural networks}
\label{subsec:training}
 
We use the Adam optimizer~\cite{kingma2017adammethodstochasticoptimization} to minimize $J^o_d$ or $J^p_d$ over the combined parameters $(\theta_y, \theta_w)$. We describe the procedure for the optimality formulation; the penalty case is analogous. 
 
Since the two networks are separate but their losses are coupled, the gradients of $J^o_d$ with respect to $\theta_y$ and $\theta_p$ are 
\begin{equation}
\frac{\partial J^o_d}{\partial \theta_y}
  \;=\; \frac{\partial R_{\mathrm{state}}}{\partial \theta_y} + \frac{\partial R_{\mathrm{adj}}}{\partial \theta_y} + \frac{\partial B_y}{\partial \theta_y},
\qquad
\frac{\partial J^o_d}{\partial \theta_p}
  \;=\; \frac{\partial R_{\mathrm{state}}}{\partial \theta_p} + \frac{\partial R_{\mathrm{adj}}}{\partial \theta_p} + \frac{\partial B_p}{\partial \theta_p}.
\label{eq:derivative-Jod}
\end{equation}
A key feature of \eqref{eq:derivative-Jod} is that $\partial J^o_d / \partial \theta_y$ depends
on the parameters $\theta_p$ of the adjoint network (through $R_{\mathrm{state}}$), and conversely $\partial J^o_d / \partial \theta_p$ depends on $\theta_y$. Even though we maintain two separate networks, their parameter updates cannot be decoupled. 

We train the two networks simultaneously rather than alternately. At each training/optimization step, the gradients in \eqref{eq:derivative-Jod} are computed for both sets of parameters by automatic differentiation in JAX  \cite{jax2018github}, and $\theta_y$ and $\theta_p$ are updated in the same Adam step. This preserves the coupling structure of the optimality system within the training/optimizer dynamics.


\subsection{Successive training for small $\eps$} 
\label{subsec:successive}

For small $\eps$, the state and adjoint layers become thin relative to the domain scale, and the residual landscape becomes strongly localized near the layers. Direct training at the target value of $\eps$ can therefore stagnate if the initial parameters and collocation set do not resolve these regions sufficiently well.

Following \cite{qiao2025two}, we use a successive training strategy. Instead of training directly at the target $\eps$ from scratch, the networks are trained along a decreasing sequence of diffusion parameters,
\[
\eps_0 > \eps_1 > \cdots > \eps_K = \eps,
\]
where the trained parameters at $\eps_k$ are used to initialize the problem at $\eps_{k+1}$. In our experiments, we typically take $\eps_0=10^{-1}$ and use eight logarithmically spaced stages from $\eps_0$ to the target $\eps$.  

At each stage $k$, we additionally apply residual-based adaptive refinement (RAR) strategy \cite{lu2021, WU2023115671} periodically during optimization. The initial collocation set contains $22500$ interior points and $1000$ boundary points. During training at a fixed $\eps_k$, RAR is performed every $20000$ Adam iterations. For the interior refinement, $20000$ candidate points are sampled in the domain, and the residual indicator
\[
\eta(x)=|R_y(x)|+|R_p(x)|
\]
is evaluated, where $R_y$ and $R_p$ denote the state and adjoint equation residuals. The $2000$ candidate points with the largest values of $\eta$ are then added to the interior collocation set.
For the boundary refinement, $20000$ candidate points are sampled over the boundary of the domain. We evaluate 
\[
\eta_{b}(x)=|R_{b,y}(x)|+|R_{b,p}(x)|,
\]
where $R_{b,y}$ and $R_{b,p}$ are the boundary residuals for the state and adjoint variables, and add the $2000$ largest residual boundary points to the boundary collocation set. 

Algorithm \ref{alg:successive} summarizes the resulting procedure. The method combines continuation in $\eps$, warm-starting from the previous stage, and adaptive refinement of collocation set. In the numerical tests, this combination stabilizes training for small diffusion parameters. It is used as a training strategy rather than as a replacement for resolving the layer.

\begin{algorithm}[t]

\caption{Successive training of double two-scale PINNs with residual-based refinement}
\label{alg:successive}
\begin{algorithmic}
\STATE \textbf{Input:} Initial interior collocation set $\mathcal{D}_{\Omega}$, initial boundary collocation set $\mathcal{D}_{\partial\Omega}$, target diffusion parameter $\varepsilon$, initial diffusion parameter $\varepsilon_0$, number of stages $K+1$, number of Adam iterations $N_k$ at stage $k$, refinement interval $N_{\rm RAR}$, number of candidate points $M$, number of added points $m$
\STATE \textbf{Output:} Trained network parameters $\theta_y,\theta_p$

\vspace{0.2cm}
\STATE \textbf{Notation:} $\mathcal{D}=\mathcal{D}_{\Omega}\cup\mathcal{D}_{\partial\Omega}$ denotes the full collocation set. $R_y$ and $R_p$ denote the interior residuals of the state and adjoint equations. $R_{b,y}$ and $R_{b,p}$ denote the boundary residuals of the state and adjoint variables.

\vspace{0.2cm}
\STATE \textbf{Initialization:}
\STATE Construct a decreasing sequence $\varepsilon_0>\varepsilon_1>\cdots>\varepsilon_K=\varepsilon$
\STATE Initialize neural networks $y_{\theta_y},p_{\theta_p}$ with random parameters
\STATE Construct the two-scale feature mapping associated with $\varepsilon_0$

\vspace{0.2cm}
\STATE \textbf{For} $k=0,\ldots,K$ \textbf{do}
\STATE \quad Set the diffusion parameter to $\varepsilon_k$
\STATE \quad Update the two-scale feature mapping using $\varepsilon_k$
\STATE \quad If $k>0$, initialize the current stage using the trained parameters from stage $k-1$

\vspace{0.1cm}
\STATE \quad \textbf{For} $n=1,\ldots,N_k$ \textbf{do}
\STATE \quad\quad Take one Adam optimization step for $\mathcal{J}_{d}^{\varepsilon_k}(\theta_y,\theta_p)$ using the current collocation set $\mathcal{D}$

\vspace{0.1cm}
\STATE \quad\quad \textbf{If} $n$ is a multiple of $N_{\rm RAR}$ \textbf{then}
\STATE \quad\quad\quad Sample $M$ candidate points in the domain $\Omega$
\STATE \quad\quad\quad Evaluate $\eta(x)=|R_y(x)|+|R_p(x)|$ on these interior candidates
\STATE \quad\quad\quad Add the $m$ interior candidates with largest $\eta(x)$ to $\mathcal{D}_{\Omega}$

\STATE \quad\quad\quad Sample $M$ candidate points on the boundary $\partial\Omega$
\STATE \quad\quad\quad Evaluate $\eta_b(x)=|R_{b,y}(x)|+|R_{b,p}(x)|$ on these boundary candidates
\STATE \quad\quad\quad Add the $m$ boundary candidates with largest $\eta_b(x)$ to $\mathcal{D}_{\partial\Omega}$
\STATE \quad\quad \textbf{End if}

\STATE \quad \textbf{End for}
\STATE \textbf{End for}
\end{algorithmic}
\end{algorithm}


\section{Numerical experiments}\label{sec:numer}

This section examines how the formulation choice affects two-scale PINN training for convection-dominated optimal control. The comparison uses the optimality-system loss \eqref{eq:jodocp} and the penalized loss \eqref{eq:jpdocp}, with separate networks for the state and for the adjoint or control. In all experiments, we use $epoch=80000$, $N_c=150\times150$ interior collocation points $\bm{x}_r$, and $N_b=250\times 4$ boundary collocation points $\bm{x}_b$. Interior collocation points are sampled uniformly, while boundary points are sampled either uniformly or from a beta distribution depending on the layer geometry. 
For the parameters $\alpha_y$ and $\alpha_p$ in \eqref{eq:Byp}, we take 
\begin{equation}\label{eq:alphayp}
    \alpha_y=\alpha_p=\begin{cases}
1000, &  epoch<10000, \\
5000, & 10000\le epoch<30000, \\
10000, & \mbox{otherwise}.
\end{cases}
\end{equation}
For the parameters $\alpha_1$ and $\alpha_2$ in \eqref{eq:jpdocp}, we take 
\begin{equation}\label{eq:alpha12}
    \alpha_1=\begin{cases}
100, &  epoch<10000, \\
500, & 10000\le epoch<30000, \\
1000, & \mbox{otherwise}
\end{cases}\quad
\alpha_2=\begin{cases}
1000, &  epoch<10000, \\
5000, & 10000\le epoch<30000, \\
10000, & \mbox{otherwise}.
\end{cases}
\end{equation}
All examples use three hidden layers with $100$ neurons per hidden layer. Unless stated otherwise, all parameters are identical across experiments. The computations are performed on WPI's Turing research cluster using one RTX PRO 6000 Blackwell GPU.

\begin{example}[Exponential Boundary Layer \cite{leykekhman2012local,liu2024multigrid}]\label{ex:bdlayer}
For the first test, take $\Omega=(0,1)^2$, $c=1$, and $\bm{\zeta}=[\sqrt{2}/2,\sqrt{2}/2]^t$. The exact solutions of \eqref{eq:speps1} are $y=\eta(x_1)\eta(x_2)$ and $p=\eta(1-x_1)\eta(1-x_2)$, where $\eta(\cdot)$ is defined by 
\begin{equation}
  \eta(z)=z^3-\frac{e^{\frac{z-1}{\varepsilon}}-e^{-1/\varepsilon}}{1-e^{-1/\varepsilon}}, \qquad \forall z \in \mathbb{R}.
\end{equation}
The state $y$ has boundary layers near $x_1=1$ and $x_2=1$, while the adjoint $p$ has boundary layers near $x_1=0$ and $x_2=0$ \cite{leykekhman2012local}. 
\end{example}

We first set $\eps=0.01$. Boundary points on all four sides are sampled from the beta distribution $f(x;0.5,0.5)$, which concentrates points near the corners $(0,0)$ and $(1,1)$ where the boundary layers occur. We adopt the beta distribution instead of a uniform distribution because it is known that boundary layers develop near the outflow boundaries. Therefore, a denser sampling near these corner points improves the resolution of the layers.

\subsection{Errors comparison with other architectures}

We first compare our architecture with three other architectures: (1) standard PINNs with two separate networks for the state and the control/adjoint state, (2) a two-scale shared neural network with two outputs, and (3) two-scale separate neural networks with the same shifting centers. Our proposed architecture uses two-scale separate neural networks with different shifting centers. We train each of these four models using three random seeds and report the average results. Table \ref{table:ex1errors} provides a comprehensive comparison of the four architectures. We can clearly see that our proposed model achieves the lowest errors in terms of $L^1$, $L^2$, and $L^\infty$ norms. The wall-clock training times are also reported. We also clearly see from Table \ref{table:ex1errors} that our proposed model achieves the best accuracy with training times comparable to those of the other models for both the optimality system and penalization approaches.

\begin{figure}[h]
    \caption{NN predictions vs Exact solutions for Example \ref{ex:bdlayer}}
    \label{nnvsexacty}
    \centering
    \subfloat[NN Solution for $y$ using \eqref{eq:jodocp}]{\includegraphics[height=2in]{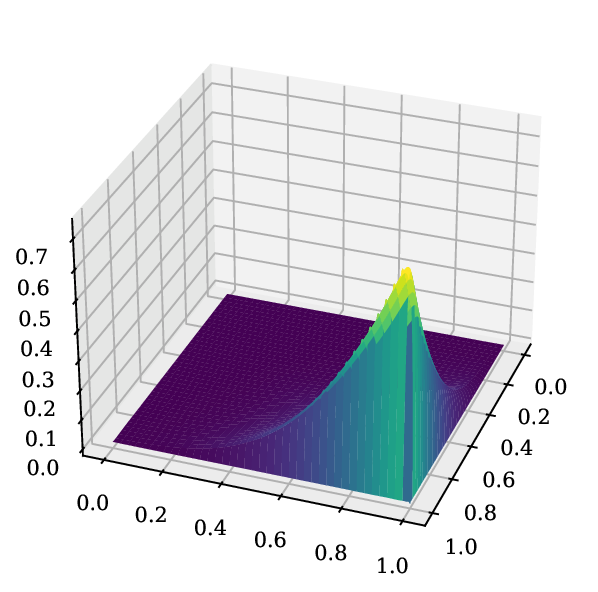}} 
    \subfloat[NN Solution for $y$ using \eqref{eq:jpdocp}]{\includegraphics[height=2in]{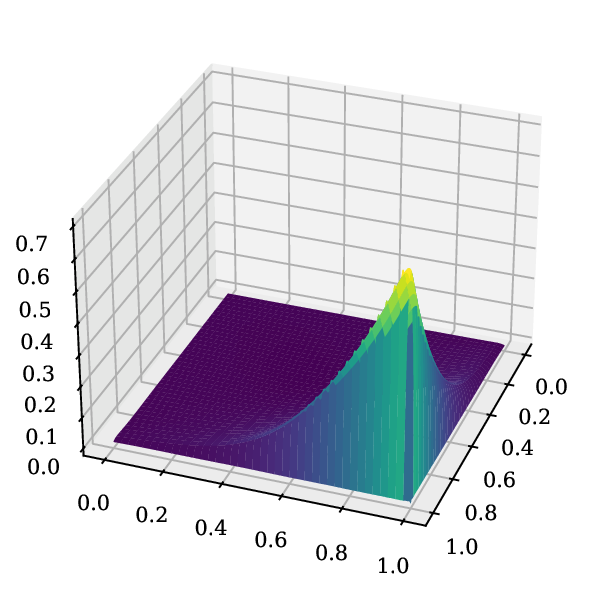}} 
    \subfloat[Exact Solution for $y$]{\includegraphics[height=2in]{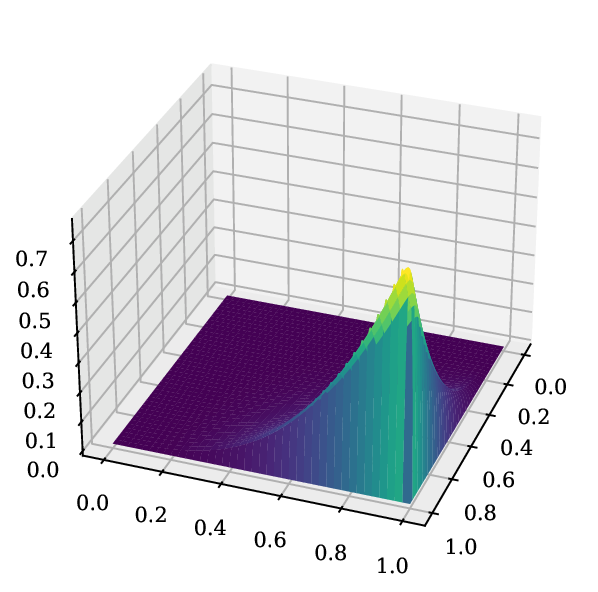}}
    \vfill
    \subfloat[NN Solution for $p$ using \eqref{eq:jodocp}]{\includegraphics[height=2in]{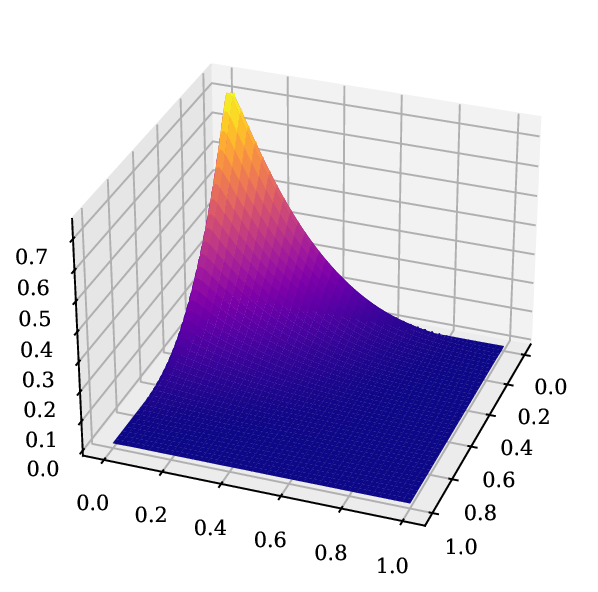}} 
    \subfloat[NN Solution for $p=-u$ using \eqref{eq:jpdocp}]{\includegraphics[height=2in]{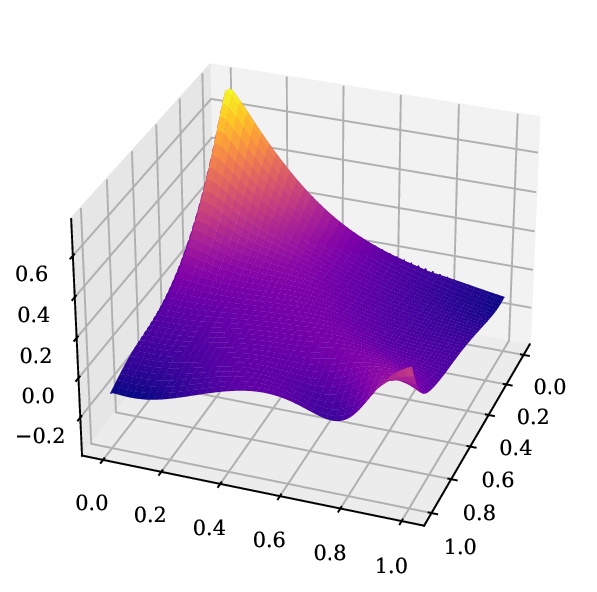}} 
    \subfloat[Exact Solution for $p$]{\includegraphics[height=2in]{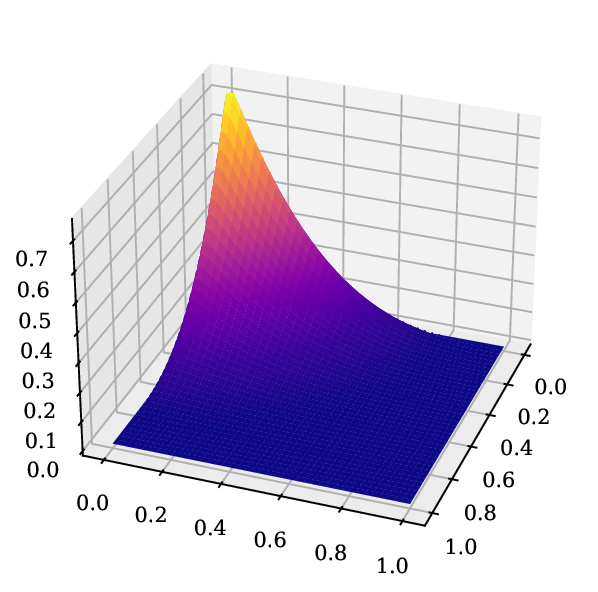}}
\end{figure}

Figure \ref{nnvsexacty} show the corresponding solution profiles. 
Both formulations capture the state boundary layer without visible spurious oscillations. The penalized loss \eqref{eq:jpdocp}, however, is less accurate for the control/adjoint variable near the boundary, whereas the optimality loss \eqref{eq:jodocp} gives smaller errors. This difference is consistent with the structure of the losses: \eqref{eq:jodocp} enforces both the state and adjoint equations, including boundary conditions for both $y$ and $p$, while \eqref{eq:jpdocp} enforces the state boundary condition directly and relies on penalty weights to recover the control behavior. 

\begin{table}[htbp]
\centering
\small
\caption{Errors and wall-clock time (mean over 3 seeds) for Example \ref{ex:bdlayer}.}\label{table:ex1errors}
\begin{tabular}{lccccccc}
\toprule
 & \multicolumn{3}{c}{State $y$} & \multicolumn{3}{c}{Control $u$} & \\
\cmidrule(lr){2-4}\cmidrule(lr){5-7}
Architecture & $L^1$ & $L^2$ & $L^\infty$ & $L^1$ & $L^2$ & $L^\infty$ & Time (s) \\
\midrule
\multicolumn{8}{l}{\textit{Optimality system}} \\
Standard PINN (separate) & 1.87e-01 & 2.77e-01 & 9.48e-01 & 1.25e-01 & 1.90e-01 & 7.46e-01 & 407.7 \\
Two-scale, shared 2-output & 1.29e-01 & 2.00e-01 & 7.84e-01 & 3.79e-02 & 5.80e-02 & 2.56e-01 & 525.4 \\
Two-scale, separate, same center & 1.50e-02 & 2.18e-02 & 8.02e-02 & 1.04e-02 & 1.56e-02 & 5.59e-02 & 401.5 \\
Two-scale, separate, different centers & 1.39e-03 & 1.86e-03 & 6.59e-03 & 1.58e-03 & 2.02e-03 & 7.99e-03 & 404.6 \\
\midrule
\multicolumn{8}{l}{\textit{Penalization}} \\
Standard PINN (separate) & 1.44e-01 & 2.21e-01 & 8.95e-01 & 5.92e-01 & 1.12e+00 & 1.24e+01 & 213.9 \\
Two-scale, shared 2-output & 3.41e-02 & 5.23e-02 & 1.88e-01 & 1.34e-01 & 2.11e-01 & 4.16e+00 & 227.5 \\
Two-scale, separate, same center & 3.49e-02 & 5.55e-02 & 2.09e-01 & 1.31e-01 & 2.60e-01 & 7.34e+00 & 212.6 \\
Two-scale, separate, different centers & 7.27e-03 & 9.83e-03 & 2.78e-02 & 2.97e-02 & 3.55e-02 & 2.18e-01 & 211.5 \\
\bottomrule
\end{tabular}
\end{table}

\subsection{Sensitivity study}

We now include a study of some important parameters involved in our numerical experiments, and see how sensitive they are in terms of errors. The tests are for the following parameters:

\begin{itemize}
\item $\varepsilon$: $\{10^{-1}, 5\times 10^{-2}, 2\times 10^{-2}, 10^{-2}, 5\times 10^{-3}, 2\times 10^{-3}, 10^{-3}, 5\times 10^{-4}\}$;
\item $\gamma$: $\{-0.5, -0.75, -1.0\}$;
\item $\beta$: $\{0.1, 1, 10\}$, using a manufactured family with the forcing recomputed consistently;
\item center placement: $x_c^y=(1-\delta,1-\delta)$, $x_c^p=(\delta,\delta)$, $\delta\in\{0,0.1,0.25,0.5\}$;
\item penalty weights: independent multipliers $s_1,s_2\in\{10,100,1000\}$ applied to the current $\alpha_y, \alpha_p$, $\alpha_1$ and $\alpha_2$ schedules.
\end{itemize}
The results are presented in Tables \ref{table:seps}-\ref{table:salpha}. For the experiments involving $\varepsilon$, neither continuation nor RAR is used. As shown in Table \ref{table:seps}, the errors behave as expected: without any specialized treatment, it is difficult to capture the boundary layer for very small values of $\varepsilon$. We also observe that the optimality system approach outperforms the penalization approach, even without continuation or RAR.
Regarding the scaling component $\gamma$, Table \ref{table:sgamma} shows that $\gamma=-0.5$ and $\gamma=-0.75$ yield similar errors, while $\gamma=-1$ produces slightly larger errors. Nevertheless, the errors remain within the same order of magnitude, indicating that our algorithm is relatively insensitive to the choice of $\gamma$.
In Table \ref{table:sbeta}, we observe that the errors decrease as $\beta$ increases. This is consistent with the well-known fact in optimal control that smaller values of $\beta$ lead to more challenging problems. Our numerical results reflect this behavior.
For the shifting centers $\delta$, Table \ref{table:sdelta} shows that our original setting has the best results. This is expected because the choice of shifting centers is closely related to the boundary layer's position, and our original setting is designed to capture this feature effectively. Nevertheless, even when the same center is used at the midpoint of the domain, our architecture can still capture the boundary layer with errors on the order of $10^{-2}$. This performance is clearly better than that of standard PINNs, as shown in Table \ref{table:ex1errors}.
Finally, for the penalty weights $\alpha_y$, $\alpha_p$ in \eqref{eq:Byp} and $\alpha_1$, $\alpha_2$ in \eqref{eq:jpdocp}, we increase all the weights simultaneously according to the prescribed schedule numbers. The base case in Table \ref{table:salpha} is the choice in \eqref{eq:alphayp} and \eqref{eq:alpha12}. As shown in Table \ref{table:salpha}, larger penalty parameters may lead to slightly worse results. Overall, however, the performance is relatively insensitive to the choice of penalty parameters.

\begin{table}[htbp]
\centering
\small
\caption{Sensitivity to $\varepsilon$.}\label{table:seps}
\begin{tabular}{lccccccc}
\toprule
 & \multicolumn{3}{c}{State $y$} & \multicolumn{3}{c}{Control $u$} & \\
\cmidrule(lr){2-4}\cmidrule(lr){5-7}
$\varepsilon$ & $L^1$ & $L^2$ & $L^\infty$ & $L^1$ & $L^2$ & $L^\infty$ & Time (s) \\
\midrule
\multicolumn{8}{l}{\textit{Optimality system}} \\
$10^{-1}$ & 3.80e-04 & 4.89e-04 & 1.46e-03 & 5.14e-04 & 7.15e-04 & 2.06e-03 & 434.5 \\
$5\times10^{-2}$ & 5.41e-04 & 7.33e-04 & 2.34e-03 & 7.10e-04 & 8.64e-04 & 2.33e-03 & 415.1 \\
$2\times10^{-2}$ & 7.97e-04 & 1.10e-03 & 4.95e-03 & 1.25e-03 & 1.54e-03 & 4.90e-03 & 415.5 \\
$10^{-2}$ & 1.15e-03 & 1.43e-03 & 5.11e-03 & 2.35e-03 & 2.83e-03 & 1.14e-02 & 414.9 \\
$5\times10^{-3}$ & 1.47e-03 & 1.82e-03 & 6.87e-03 & 1.45e-03 & 2.05e-03 & 2.27e-02 & 415.3 \\
$2\times10^{-3}$ & 1.52e-03 & 2.20e-03 & 3.11e-02 & 3.15e-03 & 4.66e-03 & 2.21e-02 & 414.6 \\
$10^{-3}$ & 7.24e-03 & 1.03e-02 & 1.08e-01 & 7.91e-03 & 1.16e-02 & 4.30e-02 & 415.2 \\
$5\times10^{-4}$ & 4.24e-03 & 6.30e-03 & 3.37e-01 & 1.66e-02 & 2.52e-02 & 1.08e-01 & 414.9 \\
\midrule
\multicolumn{8}{l}{\textit{Penalization}} \\
$10^{-1}$ & 1.84e-03 & 2.87e-03 & 1.01e-02 & 2.91e-02 & 5.72e-02 & 9.01e-01 & 223.4 \\
$5\times10^{-2}$ & 2.27e-03 & 3.34e-03 & 9.30e-03 & 2.32e-02 & 3.43e-02 & 2.66e-01 & 221.1 \\
$2\times10^{-2}$ & 4.49e-03 & 5.34e-03 & 1.12e-02 & 1.99e-02 & 2.49e-02 & 1.22e-01 & 221.2 \\
$10^{-2}$ & 7.96e-03 & 1.08e-02 & 3.41e-02 & 3.00e-02 & 3.61e-02 & 2.31e-01 & 221.3 \\
$5\times10^{-3}$ & 2.86e-02 & 4.71e-02 & 1.86e-01 & 9.86e-02 & 1.42e-01 & 5.57e-01 & 221.2 \\
$2\times10^{-3}$ & 1.38e-01 & 2.03e-01 & 7.50e-01 & 4.65e-01 & 6.07e-01 & 2.10e+00 & 220.9 \\
$10^{-3}$ & 1.75e-01 & 2.63e-01 & 1.11e+00 & 6.29e-01 & 8.79e-01 & 4.38e+00 & 221.0 \\
$5\times10^{-4}$ & 1.78e-01 & 2.85e-01 & 1.45e+00 & 7.84e-01 & 1.22e+00 & 5.70e+00 & 221.1 \\
\bottomrule
\end{tabular}
\end{table}

\begin{table}[htbp]
\centering
\small
\caption{Sensitivity to $\gamma$.}\label{table:sgamma}
\begin{tabular}{lccccccc}
\toprule
 & \multicolumn{3}{c}{State $y$} & \multicolumn{3}{c}{Control $u$} & \\
\cmidrule(lr){2-4}\cmidrule(lr){5-7}
$\gamma$ & $L^1$ & $L^2$ & $L^\infty$ & $L^1$ & $L^2$ & $L^\infty$ & Time (s) \\
\midrule
\multicolumn{8}{l}{\textit{Optimality system}} \\
$-0.5$ & 1.15e-03 & 1.43e-03 & 5.11e-03 & 2.35e-03 & 2.83e-03 & 1.14e-02 & 415.5 \\
$-0.75$ & 2.06e-03 & 2.25e-03 & 5.47e-03 & 9.62e-04 & 1.39e-03 & 5.70e-03 & 415.3 \\
$-1.0$ & 3.43e-03 & 5.19e-03 & 1.62e-02 & 3.10e-03 & 4.83e-03 & 1.75e-02 & 414.9 \\
\midrule
\multicolumn{8}{l}{\textit{Penalization}} \\
$-0.5$ & 7.96e-03 & 1.08e-02 & 3.41e-02 & 3.00e-02 & 3.61e-02 & 2.31e-01 & 221.3 \\
$-0.75$ & 7.56e-03 & 1.05e-02 & 3.02e-02 & 2.74e-02 & 3.36e-02 & 2.28e-01 & 221.4 \\
$-1.0$ & 2.41e-02 & 3.86e-02 & 1.36e-01 & 1.05e-01 & 1.44e-01 & 4.16e-01 & 221.1 \\
\bottomrule
\end{tabular}
\end{table}

\begin{table}[htbp]
\centering
\small
\caption{Sensitivity to $\beta$.}\label{table:sbeta}
\begin{tabular}{lccccccc}
\toprule
 & \multicolumn{3}{c}{State $y$} & \multicolumn{3}{c}{Control $u$} & \\
\cmidrule(lr){2-4}\cmidrule(lr){5-7}
$\beta$ & $L^1$ & $L^2$ & $L^\infty$ & $L^1$ & $L^2$ & $L^\infty$ & Time (s) \\
\midrule
\multicolumn{8}{l}{\textit{Optimality system}} \\
$0.1$ & 5.19e-03 & 7.21e-03 & 2.00e-02 & 1.52e-02 & 1.95e-02 & 9.21e-02 & 415.2 \\
$1$ & 1.15e-03 & 1.43e-03 & 5.11e-03 & 2.35e-03 & 2.83e-03 & 1.14e-02 & 415.1 \\
$10$ & 1.67e-03 & 2.13e-03 & 7.17e-03 & 2.34e-04 & 2.95e-04 & 7.88e-04 & 415.1 \\
\midrule
\multicolumn{8}{l}{\textit{Penalization}} \\
$0.1$ & 4.27e-02 & 8.86e-02 & 5.03e-01 & 2.67e-01 & 5.55e-01 & 3.91e+00 & 221.2 \\
$1$ & 7.96e-03 & 1.08e-02 & 3.41e-02 & 3.00e-02 & 3.61e-02 & 2.31e-01 & 221.3 \\
$10$ & 1.22e-03 & 1.51e-03 & 6.70e-03 & 6.71e-03 & 9.76e-03 & 6.81e-02 & 221.3 \\
\bottomrule
\end{tabular}
\end{table}

\begin{table}[htbp]
\centering
\small
\caption{Sensitivity to center offset $\delta$.}\label{table:sdelta}
\begin{tabular}{lccccccc}
\toprule
 & \multicolumn{3}{c}{State $y$} & \multicolumn{3}{c}{Control $u$} & \\
\cmidrule(lr){2-4}\cmidrule(lr){5-7}
$\delta$ & $L^1$ & $L^2$ & $L^\infty$ & $L^1$ & $L^2$ & $L^\infty$ & Time (s) \\
\midrule
\multicolumn{8}{l}{\textit{Optimality system}} \\
$0$ & 1.15e-03 & 1.43e-03 & 5.11e-03 & 2.35e-03 & 2.83e-03 & 1.14e-02 & 415.4 \\
$0.1$ & 3.26e-03 & 4.71e-03 & 1.77e-02 & 2.44e-03 & 3.00e-03 & 8.87e-03 & 415.0 \\
$0.25$ & 9.00e-03 & 1.28e-02 & 4.41e-02 & 5.05e-03 & 6.97e-03 & 2.19e-02 & 415.5 \\
$0.5$ & 1.22e-02 & 1.73e-02 & 5.63e-02 & 8.34e-03 & 1.17e-02 & 3.89e-02 & 420.8 \\
\midrule
\multicolumn{8}{l}{\textit{Penalization}} \\
$0$ & 7.96e-03 & 1.08e-02 & 3.41e-02 & 3.00e-02 & 3.61e-02 & 2.31e-01 & 221.3 \\
$0.1$ & 1.07e-02 & 1.52e-02 & 5.29e-02 & 4.18e-02 & 5.52e-02 & 5.48e-01 & 221.5 \\
$0.25$ & 1.24e-02 & 1.89e-02 & 6.28e-02 & 5.03e-02 & 6.48e-02 & 2.85e-01 & 221.4 \\
$0.5$ & 3.02e-02 & 4.91e-02 & 1.89e-01 & 1.30e-01 & 2.89e-01 & 9.51e+00 & 222.9 \\
\bottomrule
\end{tabular}
\end{table}

\begin{table}[htbp]
\centering
\small
\caption{Sensitivity to penalty-weight increment on $\alpha_y$ $\alpha_p$, $\alpha_1,\alpha_2$.}\label{table:salpha}
\begin{tabular}{lccccccc}
\toprule
 & \multicolumn{3}{c}{State $y$} & \multicolumn{3}{c}{Control $u$} & \\
\cmidrule(lr){2-4}\cmidrule(lr){5-7}
increment & $L^1$ & $L^2$ & $L^\infty$ & $L^1$ & $L^2$ & $L^\infty$ & Time (s) \\
\midrule
\multicolumn{8}{l}{\textit{Optimality system}} \\
base & 8.52e-04 & 1.23e-03 & 5.08e-03 & 8.65e-04 & 1.20e-03 & 6.31e-03 & 445.0 \\
$+10$ & 1.78e-03 & 2.26e-03 & 6.86e-03 & 2.52e-03 & 2.71e-03 & 7.69e-03 & 425.2 \\
$+100$ & 1.57e-03 & 1.84e-03 & 5.90e-03 & 1.39e-03 & 2.06e-03 & 8.67e-03 & 425.2 \\
$+1000$ & 1.84e-03 & 2.35e-03 & 7.99e-03 & 2.60e-03 & 3.93e-03 & 1.57e-02 & 425.4 \\
\midrule
\multicolumn{8}{l}{\textit{Penalization}} \\
base & 6.99e-03 & 8.90e-03 & 2.16e-02 & 3.18e-02 & 3.82e-02 & 2.34e-01 & 226.3 \\
$+10$ & 9.80e-03 & 1.30e-02 & 3.77e-02 & 3.99e-02 & 4.81e-02 & 2.37e-01 & 223.9 \\
$+100$ & 1.03e-02 & 1.39e-02 & 5.13e-02 & 4.08e-02 & 5.21e-02 & 2.64e-01 & 224.1 \\
$+1000$ & 2.83e-02 & 4.46e-02 & 1.63e-01 & 1.08e-01 & 1.50e-01 & 9.24e-01 & 223.9 \\
\bottomrule
\end{tabular}
\end{table}

\subsection{Continuation and RAR study}

\begin{figure}
 \caption{NN predictions vs Exact solutions using \eqref{eq:jodocp} for Example \ref{ex:bdlayer} with $\eps=5\times10^{-4}$}
    \label{nnvsexactysmall}
    \centering
    \subfloat[NN Solution for $y$]{\includegraphics[height=2in]{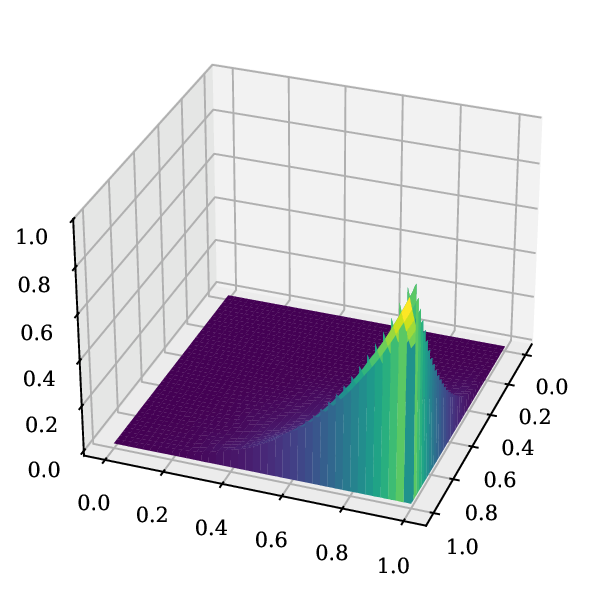}} 
    \subfloat[NN Solution for $p$]{\includegraphics[height=2in]{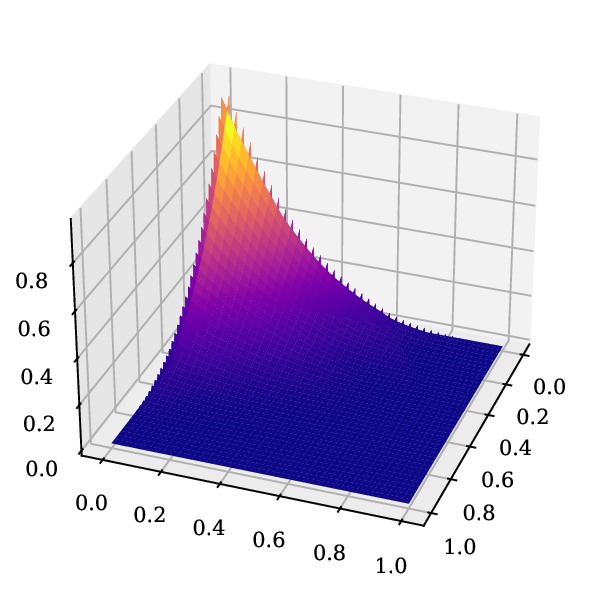}} 
\end{figure}

We next set $\eps = 5\times10^{-4}$ and apply the successive training strategy in Algorithm \ref{alg:successive}. Figure \ref{nnvsexactysmall} shows that the network predictions resolve the exponential boundary layers without visible spurious oscillations. The absolute errors for both $y$ and $p$ are of order $10^{-3}$. This experiment supports the use of two-scale continuation as a practical training strategy for small diffusion, while leaving open the question of how best to tune continuation and refinement for still smaller values of $\eps$.
We also investigate the effects of continuation and RAR in Table \ref{table:conRAR}. We observe that continuation alone leads to a noticeable improvement, while combining continuation with RAR gives slightly better overall results, although the additional improvement is not significant. 

\begin{table}[htbp]
\centering
\small
\caption{Optimality system, $\varepsilon=5\times10^{-4}$: effect of continuation and RAR.}\label{table:conRAR}
\begin{tabular}{lccccccc}
\toprule
 & \multicolumn{3}{c}{State $y$} & \multicolumn{3}{c}{Control $u$} & \\
\cmidrule(lr){2-4}\cmidrule(lr){5-7}
Strategy & $L^1$ & $L^2$ & $L^\infty$ & $L^1$ & $L^2$ & $L^\infty$ & Time (s) \\
\midrule
No continuation, no RAR & 8.23e-02 & 1.75e-01 & 1.31e+00 & 6.29e-02 & 1.40e-01 & 9.70e-01 & 557.8 \\
Continuation, no RAR & 1.29e-03 & 2.56e-03 & 1.07e-01 & 1.84e-03 & 2.57e-03 & 2.84e-02 & 2878.8 \\
Continuation + RAR & 2.35e-03 & 3.58e-03 & 1.28e-02 & 3.88e-03 & 6.03e-03 & 2.08e-02 & 3724.7 \\
\bottomrule
\end{tabular}
\end{table}

\begin{example}[Interior Layer \cite{leykekhman2012local,liu2024multigrid}]\label{ex:inlayer}
For the second test, take $\Omega=(0,1)^2$, $c=1$, $\bm{\zeta}=[1, 0]^t$, and $\eps=0.01$. The exact solution of \eqref{eq:speps1} is 
    \begin{equation}
        y=(1-x_1)^3\arctan(\frac{x_2-0.5}{\eps})\quad\text{and}\quad p=x_1(1-x_1)x_2(1-x_2).
    \end{equation}
    The exact state $y$ has an interior layer along the line $x_2=0.5$ for small $\eps$. 
    
\end{example}

\begin{figure}[h]
    \caption{NN predictions vs Exact solutions for Example \ref{ex:inlayer}}
    \label{nnvsexacty2}
    \centering
    \subfloat[NN Solution for $y$ using using \eqref{eq:jodocp}]{\includegraphics[height=2in]{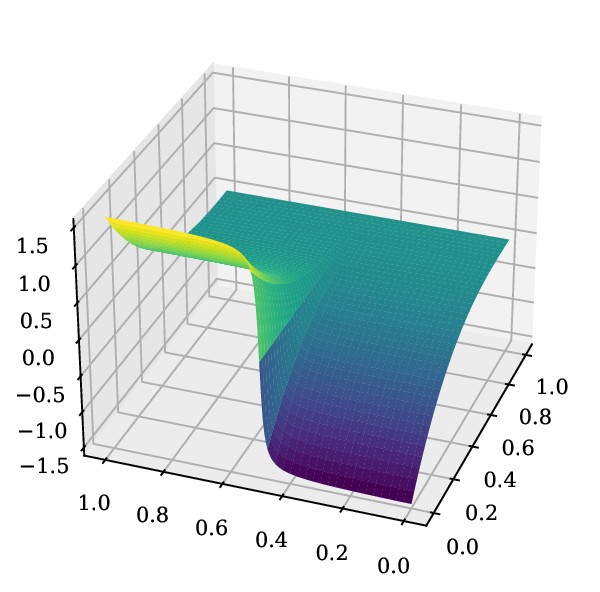}} 
    \subfloat[NN Solution for $y$ using \eqref{eq:jpdocp}]{\includegraphics[height=2in]{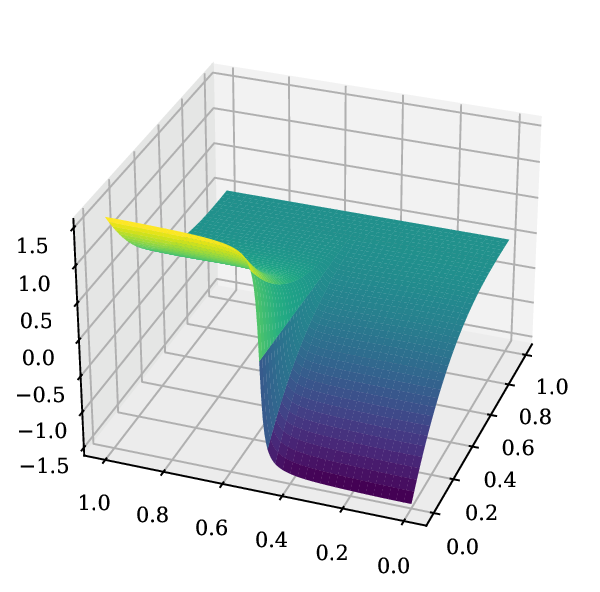}} 
    \subfloat[Exact Solution for $y$]{\includegraphics[height=2in]{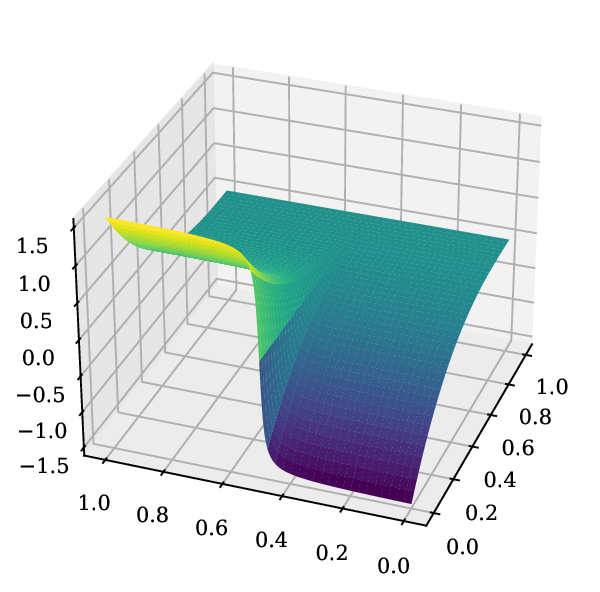}}
    \vfill
    \subfloat[NN Solution for $p=-u$ using \eqref{eq:jodocp}]{\includegraphics[height=2in]{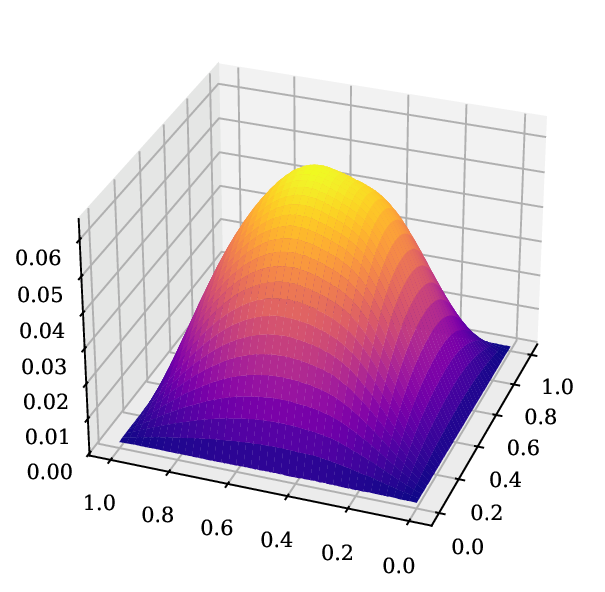}} 
    \subfloat[NN Solution for $p=-u$ using \eqref{eq:jpdocp}]{\includegraphics[height=2in]{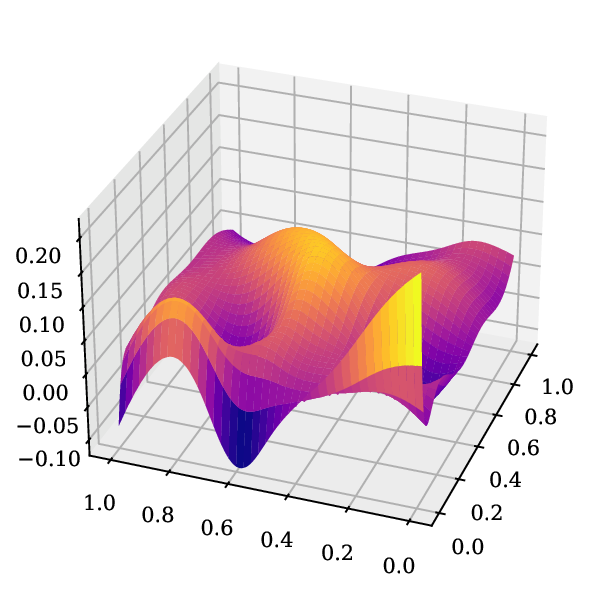}} 
    \subfloat[Exact Solution for $p=-u$]{\includegraphics[height=2in]{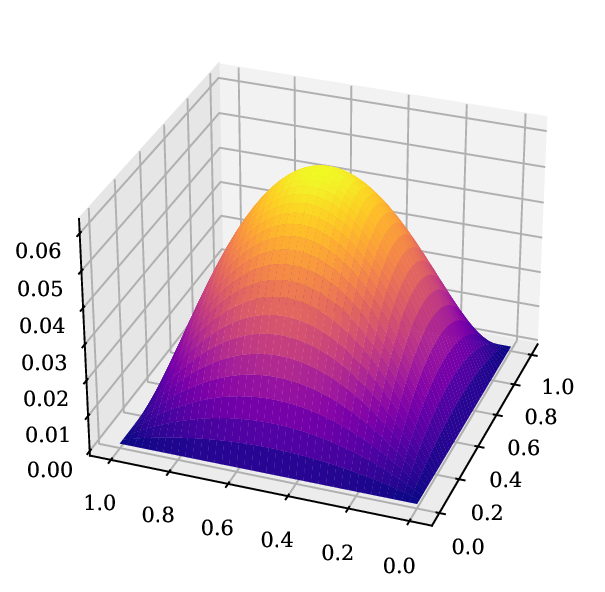}}
\end{figure}

Boundary points are sampled uniformly along each side of the square. Figure \ref{nnvsexacty2} shows the solution profiles. With the optimality loss \eqref{eq:jodocp}, the two-scale architecture resolves the interior layer in $y$ and accurately approximates the smooth adjoint solution $p$, achieving errors on the order of $10^{-3}$. With the penalized loss \eqref{eq:jpdocp}, the approximation of the control variable $u$ fails to converge in this test, consistent with the results for the boundary-layer example.

Table \ref{table:ex2errors} presents a comprehensive comparison of models similar to those considered in Example \ref{ex:bdlayer}. We observe similar convergence behaviors, and our proposed architecture continues to achieve the best performance. Similar cost-accuracy trends are also observed in Table \ref{table:ex2errors}.

\begin{table}[htbp]
\centering
\small
\caption{Errors and wall-clock time (mean over 3 seeds) for Example \ref{ex:inlayer}.}\label{table:ex2errors}
\begin{tabular}{lccccccc}
\toprule
 & \multicolumn{3}{c}{State $y$} & \multicolumn{3}{c}{Control $u$} & \\
\cmidrule(lr){2-4}\cmidrule(lr){5-7}
Architecture & $L^1$ & $L^2$ & $L^\infty$ & $L^1$ & $L^2$ & $L^\infty$ & Time (s) \\
\midrule
\multicolumn{8}{l}{\textit{Optimality system}} \\
Standard PINN (separate) & 1.74e-03 & 2.33e-03 & 1.80e-02 & 6.89e-03 & 8.15e-03 & 1.55e-02 & 406.8 \\
Two-scale, shared 2-output & 1.45e-03 & 1.85e-03 & 6.49e-03 & 2.40e-03 & 2.99e-03 & 7.58e-03 & 524.2 \\
Two-scale, separate, same center & 1.23e-03 & 1.50e-03 & 4.19e-03 & 7.98e-04 & 1.01e-03 & 2.90e-03 & 401.7 \\
\midrule
\multicolumn{8}{l}{\textit{Penalization}} \\
Standard PINN (separate) & 5.20e-03 & 6.83e-03 & 1.79e-02 & 2.89e-02 & 3.73e-02 & 1.01e-01 & 213.7 \\
Two-scale, shared 2-output & 1.92e-03 & 2.34e-03 & 7.24e-03 & 1.16e-02 & 1.80e-02 & 1.39e-01 & 227.2 \\
Two-scale, separate, same center & 2.28e-03 & 2.91e-03 & 7.96e-03 & 1.26e-02 & 1.62e-02 & 1.04e-01 & 212.7 \\
\bottomrule
\end{tabular}
\end{table}

\begin{example}[Parabolic Boundary Layers \cite{roos2008robust}] \label{ex:paralayer}
For the final test, take $\Omega = (0,1)^2$, $c = 1$, $\bm{\zeta} = [1, 0]^t$, and $\eps=0.01$, with $y_d=1$ and $f=1$. No closed-form solution is available. The state $y$ has an exponential boundary layer at the outflow boundary $x_1=1$ and parabolic boundary layers at the characteristic boundaries $x_2=0$ and $x_2=1$, while the adjoint state $p$ and control $u$ have an exponential boundary layer at $x_1=0$.
\end{example}

\begin{figure}[h]
    \caption{Comparison between NN predictions and EAFE solutions using \eqref{eq:jodocp} for Example \ref{ex:paralayer}}
    \label{nnpredict3}
    \centering
    \subfloat[NN Solution for $y$ using \eqref{eq:jodocp}]{\includegraphics[height=2in]{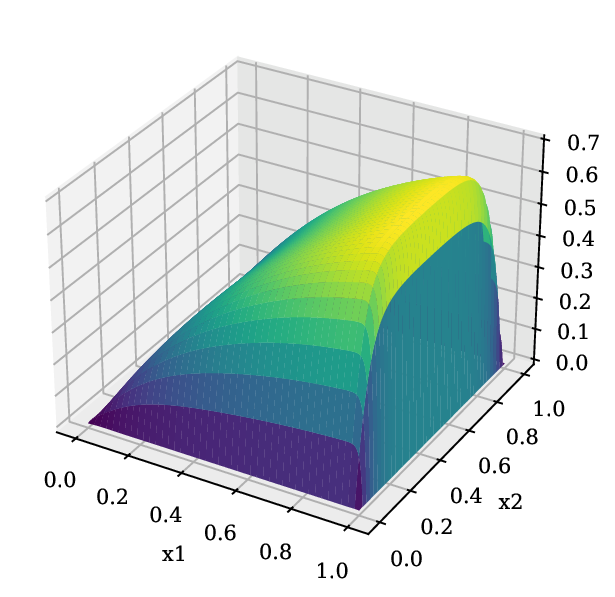}}
    \subfloat[EAFE Solution for $y$]{\includegraphics[height=2in]{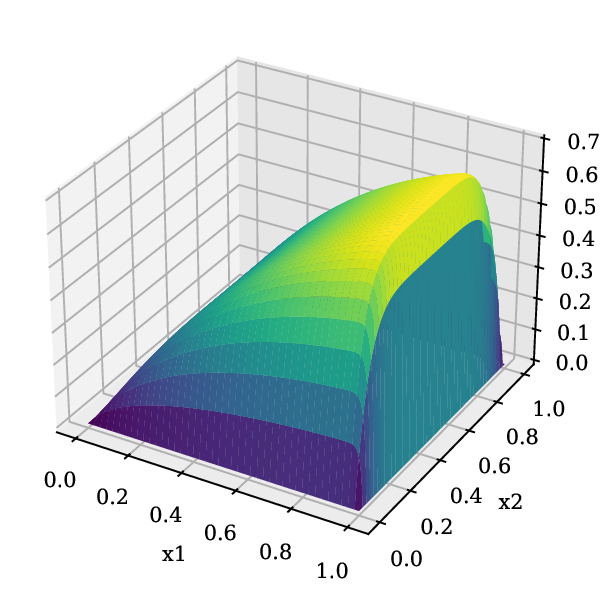}}
    \vfill
    \subfloat[NN Solution for $p$ using \eqref{eq:jodocp}]
    {\includegraphics[height=2in]{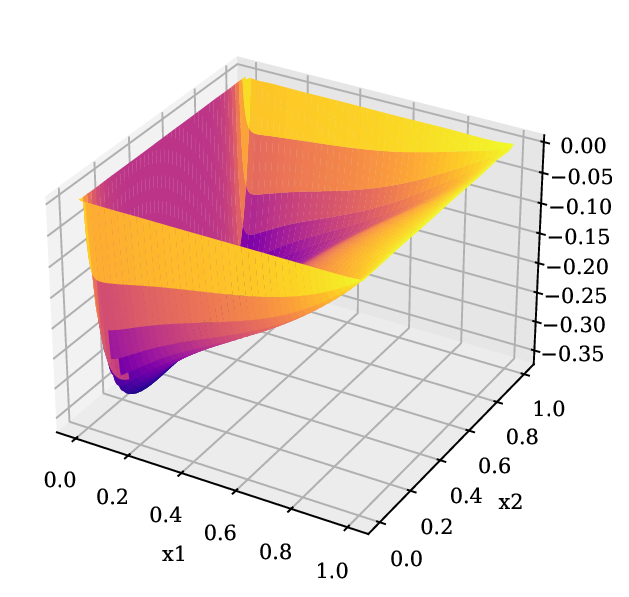}}
    \subfloat[EAFE Solution for $p$]{\includegraphics[height=2in]{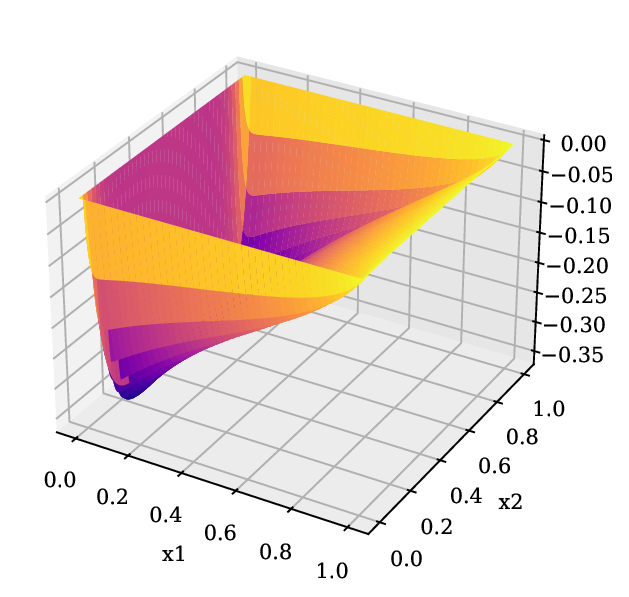}}
\end{figure}

Boundary points are sampled uniformly along all four sides. Since no exact solution is available, reference solutions are computed on a structured mesh with mesh size $h=2^{-7}$ using the edge-average finite element (EAFE) method \cite{xu1999monotone,jeong2025monotone}. EAFE stabilizes convection-dominated discretizations, captures boundary layers, and offers theoretical convergence guarantees.
Figure \ref{nnpredict3} compares the neural network predictions with the corresponding EAFE solutions, while Table \ref{table:ex3errors} presents a comprehensive comparison of different architectures. The qualitative behavior is consistent with that observed in the first two examples: the two-scale representation captures the layer structure, while the choice of formulation affects the accuracy of the adjoint/control approximation. Our proposed architectures achieve the best overall performance, and the cost-accuracy trends are also similar to those observed in the first two examples.

\begin{table}[htbp]
\centering
\small
\caption{Errors against the FEM solution and wall-clock time (mean over 3 seeds) for Example \ref{ex:paralayer}.}\label{table:ex3errors}
\begin{tabular}{lccccccc}
\toprule
 & \multicolumn{3}{c}{State $y$} & \multicolumn{3}{c}{Control $u$} & \\
\cmidrule(lr){2-4}\cmidrule(lr){5-7}
Architecture & $L^1$ & $L^2$ & $L^\infty$ & $L^1$ & $L^2$ & $L^\infty$ & Time (s) \\
\midrule
\multicolumn{8}{l}{\textit{Optimality system}} \\
Standard PINN (separate) & 2.52e-01 & 3.11e-01 & 6.57e-01 & 1.58e-01 & 1.88e-01 & 3.64e-01 & 409.9 \\
Two-scale, shared 2-output & 1.74e-01 & 2.17e-01 & 4.58e-01 & 1.37e-01 & 1.71e-01 & 3.57e-01 & 530.1 \\
Two-scale, separate, same center & 1.60e-02 & 2.04e-02 & 5.36e-02 & 1.38e-01 & 1.72e-01 & 3.55e-01 & 404.2 \\
Two-scale, separate, different centers & 1.93e-03 & 2.50e-03 & 1.07e-02 & 4.45e-03 & 5.83e-03 & 2.51e-02 & 407.8 \\
\midrule
\multicolumn{8}{l}{\textit{Penalization}} \\
Standard PINN (separate) & 2.90e-01 & 3.45e-01 & 6.65e-01 & 1.05e+00 & 1.09e+00 & 1.70e+00 & 215.9 \\
Two-scale, shared 2-output & 2.50e-01 & 3.09e-01 & 6.54e-01 & 9.55e-01 & 1.04e+00 & 2.10e+00 & 229.6 \\
Two-scale, separate, same center & 2.48e-01 & 3.09e-01 & 6.51e-01 & 9.42e-01 & 1.03e+00 & 2.31e+00 & 214.6 \\
Two-scale, separate, different centers & 2.47e-01 & 3.07e-01 & 6.47e-01 & 9.33e-01 & 1.03e+00 & 2.41e+00 & 213.8 \\
\bottomrule
\end{tabular}
\end{table}

\section{Concluding remarks}\label{sec:summary}

This work studies PINN formulation choice and two-scale continuation for elliptic distributed optimal control problems constrained by convection-dominated equations. The method extends the two-scale PINN architecture of \cite{qiao2025two} from single equations to coupled optimality systems and penalized optimal control formulations. Separate two-scale networks are used for the state and for the adjoint or control, so that their stretched coordinates can be centered near different layer regions. Simultaneous training preserves the coupling in the optimality-system formulation, while continuation in $\eps$ provides a practical way to train in more singular regimes.

The experiments indicate that the optimality-system loss is more reliable for recovering the adjoint/control in the tested layer-dominated regimes, whereas the penalized loss is simpler but more sensitive to penalty choices. This conclusion is empirical and problem dependent. Its role is to clarify the tradeoff between enforcing the adjoint equation and relying on penalization, rather than to rule out penalized PINNs for broader classes of optimal control problems.

Several extensions remain natural. For nonlinear state equations, such as the Navier--Stokes equations, the corresponding adjoint systems are more involved; linearization offers one possible route, as discussed in \cite{casas2007error}. State and control constraints can also be incorporated. Classical numerical methods for such constrained problems have been developed in \cite{liu2024multigrid,liu2024discontinuous,brenner2021p1,brenner2023multigrid}, but pointwise state constraints remain computationally demanding. PINN formulations based on either penalization or optimality systems provide a flexible framework for investigating these problems, and their stability and accuracy in singularly perturbed regimes are topics for future work.
 
\section*{Funding statement}
There is no funding for this research.

\section*{Declaration of competing interest}
The authors declare that they have no known competing financial interests or personal relationships that could have influenced the work reported in this manuscript. This research was conducted independently, without any commercial or financial support that could be construed as a potential conflict of interest.
All authors have reviewed and approved the final version of the manuscript and agree with its submission. 
The authors affirm that the work represents their original research, has not been published previously, and is not under consideration elsewhere.
\section*{Data availability statements}
The code and data that support the findings of this study are available on request from the corresponding author.
\section*{Declaration of generative AI and AI-assisted technologies in the manuscript preparation}
During the preparation of this manuscript, the authors used generative AI-based tools solely for language polishing, proofreading, and improving the clarity and readability of the text.
The authors affirm that they take full responsibility for the content of the manuscript. All scientific ideas, methodologies, algorithms, experimental designs, benchmark tests, and conclusions presented in this work were conceived, developed, and validated by the authors. The use of AI tools did not influence the scientific content or the interpretation of the results.


\bibliographystyle{abbrv}
\bibliography{references,boundarylayer,control_pinn}

@article{gao2025more,
  title={More Consistent Accuracy {PINN} via Alternating Easy-Hard Training},
  author={Gao, Zhaoqian and Yang, Min},
  journal={arXiv:2512.17607},
  year={2025}
}

@article{cao2024multistep,
  title={Multistep asymptotic pre-training strategy based on PINNs for solving steep boundary singular perturbation problems},
  author={Cao, Fujun and Gao, Fei and Yuan, Dongfang and Liu, Junmin},
  journal={Computer Methods in Applied Mechanics and Engineering},
  volume={431},
  pages={117222},
  year={2024},
  publisher={Elsevier}
}

@article{wang2024aspinn,
  title={ASPINN: An asymptotic strategy for solving singularly perturbed differential equations},
  author={Wang, Sen and Zhao, Peizhi and Song, Tao},
  journal={arXiv:2409.13185},
  year={2024}
}

@article{cao2023physics,
  title={Physics-informed neural networks with parameter asymptotic strategy for learning singularly perturbed convection-dominated problem},
  author={Cao, Fujun and Gao, Fei and Guo, Xiaobin and Yuan, Dongfang},
  journal={Computers \& Mathematics with Applications},
  volume={150},
  pages={229--242},
  year={2023},
  publisher={Elsevier}
}

@article{wang2024less,
  title={Less emphasis on hard regions: curriculum learning of {PINN}s for singularly perturbed convection-diffusion-reaction problems},
  author={Wang, Yufeng and Xu, Cong and Yang, Min and Zhang, Jin},
  journal={East Asian Journal on Applied Mathematics},
  volume={14},
  number={1},
  pages={104--123},
  year={2024}
}

@article{munzer2022curriculum,
  title={A curriculum-training-based strategy for distributing collocation points during physics-informed neural network training},
  author={M{\"u}nzer, Marcus and Bard, Chris},
  journal={arXiv:2211.11396},
  year={2022}
  }

@article{wang2024general,
  title={General-Kindred physics-informed neural network to the solutions of singularly perturbed differential equations},
  author={Wang, Sen and Zhao, Peizhi and Ma, Qinglong and Song, Tao},
  journal={Physics of Fluids},
  volume={36},
  number={11},
  year={2024},
  publisher={AIP Publishing}
}

@article{jin2023solving,
  title={Solving Elliptic Optimal Control Problems Using Physics Informed Neural Networks},
  author={Jin, Bangti and Sau, Ramesh and Yin, Luowei and Zhou, Zhi},
  journal={arXiv:2308.11925},
  year={2023}
}

@inproceedings{Nzoyem2023comparison,
author = {Nzoyem Ngueguin, Roussel Desmond and Barton, David A.W. and Deakin, Tom},
title = {A Comparison of Mesh-Free Differentiable Programming and Data-Driven Strategies for Optimal Control under PDE Constraints},
year = {2023},
isbn = {9798400707858},
publisher = {Association for Computing Machinery},
address = {New York, NY, USA},
url = {https://doi.org/10.1145/3624062.3626078},
doi = {10.1145/3624062.3626078},
booktitle = {Proceedings of the SC '23 Workshops of the International Conference on High Performance Computing, Network, Storage, and Analysis},
pages = {21–28},
numpages = {8},
location = {Denver, CO, USA},
series = {SC-W'23}
}

@article{barry2025physics,
  title={Physics-informed neural networks for {PDE}-constrained optimization and control},
  author={Barry-Straume, Jostein and Sarshar, Arash and Popov, Andrey A and Sandu, Adrian},
  journal={Communications on Applied Mathematics and Computation},
  pages={1--24},
  year={2025},
  publisher={Springer}
}

@article{cao2025adversarial,
  title={Adversarial physics-informed neural networks with hard constraints for optimal control of {PDE}s},
  author={Cao, Yuandong and So, Chi Chiu and Dai, Yifan and Yung, Siu Pang and Wang, Jun-Min},
  journal={Journal of Computational Physics},
  pages={114307},
  year={2025},
  publisher={Elsevier}
}

@article{wang2026optimal,
  title={When Optimal Control Meets Neural Network: A Comprehensive Survey},
  author={Wang, Xinwei and Dou, Yuqi and Yi, Xueling and Zhang, Yepeng and Li, Xin and Li, Bai and Peng, Haijun and Wang, Lei and Teo, Kok Lay},
  journal={Archives of Computational Methods in Engineering},
  pages={1--56},
  year={2026},
  publisher={Springer}
}

@article{yin2024aonn,
  title={{AONN}: An adjoint-oriented neural network method for all-at-once solutions of parametric optimal control problems},
  author={Yin, Pengfei and Xiao, Guangqiang and Tang, Kejun and Yang, Chao},
  journal={SIAM Journal on Scientific Computing},
  volume={46},
  number={1},
  pages={C127--C153},
  year={2024},
  publisher={SIAM}
}

@article{wang2024aonn,
  title={{AONN}-2: An adjoint-oriented neural network method for {PDE}-constrained shape optimization},
  author={Wang, Xili and Yin, Pengfei and Zhang, Bo and Yang, Chao},
  journal={Journal of Computational Physics},
  volume={513},
  pages={113160},
  year={2024},
  publisher={Elsevier}
}

@article{dupret2026deep,
  title={Deep learning for high-dimensional continuous-time stochastic optimal control without explicit solution},
  author={Dupret, Jean-Loup and Hainaut, Donatien},
  journal={Operations Research},
  year={2026},
  publisher={INFORMS},
  url = {https://doi.org/10.1287/opre.2024.1102},
}

@article{qiao2026two,
  title={Two-scale Neural Networks for Singularly Perturbed Dynamical Systems with Multiple Parameters},
  author={Zhuang, Qiao and Wang, Taorui and Wanjiku, Rita and Bani-Yaghoub, Majid and Zhang, Zhongqiang},
  journal={arXiv:2605.02799},
  year={2026}
}

@article{brenner2023multigrid,
  title={Multigrid methods for an elliptic optimal control problem with pointwise state constraints},
  author={Brenner, S. C. and Liu, S. and Sung, L.-Y.},
  journal={Results in Applied Mathematics},
  volume={17},
  pages={100356},
  year={2023},
  publisher={Elsevier}
}

@article{brenner2021p1,
  title={A $P_1$ Finite Element Method for a Distributed Elliptic Optimal Control Problem with a General State Equation and Pointwise State Constraints},
  author={Brenner, S. C. and Liu, S. and Sung, L.-Y.},
  journal={Computational Methods in Applied Mathematics},
  volume={21},
  number={4},
  pages={777--790},
  year={2021},
  publisher={De Gruyter}
}

@article{liu2024discontinuous,
  title={Discontinuous Galerkin methods for an elliptic optimal control problem with a general state equation and pointwise state constraints},
  author={Liu, S. and Tan, Z. and Zhang, Y.},
  journal={Journal of Computational and Applied Mathematics},
  volume={437},
  pages={115494},
  year={2024},
  publisher={Elsevier}
}

@article{casas2007error,
  title={Error estimates for the numerical approximation of a distributed control problem for the steady-state Navier--Stokes equations},
  author={Casas, E. and Mateos, M. and Raymond, J.-P.},
  journal={SIAM Journal on Control and Optimization},
  volume={46},
  number={3},
  pages={952--982},
  year={2007},
  publisher={SIAM}
}

@software{jax2018github,
  author = {James Bradbury and Roy Frostig and Peter Hawkins and Matthew James Johnson and Yash Katariya and Chris Leary and Dougal Maclaurin and George Necula and Adam Paszke and Jake Vander{P}las and Skye Wanderman-{M}ilne and Qiao Zhang},
  title = {{JAX}: composable transformations of {P}ython+{N}um{P}y programs},
  url = {http://github.com/jax-ml/jax},
  version = {0.3.13},
  year = {2018},
}

@article{lu2021,
author = {Lu, Lu and Meng, Xuhui and Mao, Zhiping and Karniadakis, George Em},
title = {DeepXDE: A Deep Learning Library for Solving Differential Equations},
journal = {SIAM Review},
volume = {63},
number = {1},
pages = {208-228},
year = {2021},
doi = {10.1137/19M1274067},

URL = { 
    
        https://doi.org/10.1137/19M1274067
    
    

},
eprint = { 
    
        https://doi.org/10.1137/19M1274067
    
    

}
}

@article{WU2023115671,
title = {A comprehensive study of non-adaptive and residual-based adaptive sampling for physics-informed neural networks},
journal = {Computer Methods in Applied Mechanics and Engineering},
volume = {403},
pages = {115671},
year = {2023},
issn = {0045-7825},
doi = {https://doi.org/10.1016/j.cma.2022.115671},
url = {https://www.sciencedirect.com/science/article/pii/S0045782522006260},
author = {Chenxi Wu and Min Zhu and Qinyang Tan and Yadhu Kartha and Lu Lu}
}

@article{BROOKS1982199,
title = {Streamline upwind/Petrov-Galerkin formulations for convection dominated flows with particular emphasis on the incompressible Navier-Stokes equations},
journal = {Computer Methods in Applied Mechanics and Engineering},
volume = {32},
number = {1},
pages = {199-259},
year = {1982},
issn = {0045-7825},
doi = {https://doi.org/10.1016/0045-7825(82)90071-8},
url = {https://www.sciencedirect.com/science/article/pii/0045782582900718},
author = {Alexander N. Brooks and Thomas J.R. Hughes}
}

@inproceedings{Hughes1979MULTIDIMENSIONALUS,
  title={MULTI-DIMENSIONAL UPWIND SCHEME WITH NO CROSSWIND DIFFUSION.},
  author={Thomas Joseph Robert Hughes and Alexander Nelson Brooks},
  year={1979},
  url={https://api.semanticscholar.org/CorpusID:118669634}
}

@article{doi:10.1137/S0036142997316712,
author = {Cockburn, Bernardo and Shu, Chi-Wang},
title = {The Local Discontinuous Galerkin Method for Time-Dependent Convection-Diffusion Systems},
journal = {SIAM Journal on Numerical Analysis},
volume = {35},
number = {6},
pages = {2440-2463},
year = {1998},
doi = {10.1137/S0036142997316712},

URL = {  
        https://doi.org/10.1137/S0036142997316712
    
}, 
}

@article{doi:10.1137/S0036142900374111,
author = {Houston, Paul and Schwab, Christoph and S\"{u}li, Endre},
title = {Discontinuous hp-Finite Element Methods for Advection-Diffusion-Reaction Problems},
journal = {SIAM Journal on Numerical Analysis},
volume = {39},
number = {6},
pages = {2133-2163},
year = {2002},
doi = {10.1137/S0036142900374111},

URL = { 
    
        https://doi.org/10.1137/S0036142900374111
    
    

},
eprint = { 
    
        https://doi.org/10.1137/S0036142900374111
    
    

}
}

@article{chen2018hdg,
  title={An {HDG} method for distributed control of convection diffusion {PDE}s},
  author={Chen, G. and Hu, W. and Shen, J. and Singler, J. R. and Zhang, Y. and Zheng, X.},
  journal={Journal of Computational and Applied Mathematics},
  volume={343},
  pages={643--661},
  year={2018},
  publisher={Elsevier}
}

@article{liu2025balancing,
  title={A balancing domain decomposition by constraints preconditioner for a hybridizable discontinuous {G}alerkin discretization of an elliptic optimal control problem},
  author={Liu, S. and Zhang, J.},  journal={arXiv:2504.02072},
  year={2025}
}

@article{liu2025convergence,
  title={Convergence analysis of a balancing domain decomposition method for an elliptic optimal control problem with {HDG} discretizations},
  author={Liu, S. and Zhang, J.},
  journal={ESAIM: Mathematical Modelling and Numerical Analysis},
  year={2026}
}

@inproceedings{nitsche1971variationsprinzip,
  title={{\"U}ber ein Variationsprinzip zur L{\"o}sung von Dirichlet-Problemen bei Verwendung von Teilr{\"a}umen, die keinen Randbedingungen unterworfen sind},
  author={Nitsche, J.},
  booktitle={Abhandlungen aus dem mathematischen Seminar der Universit{\"a}t Hamburg},
  volume={36},
  number={1},
  pages={9--15},
  year={1971},
  organization={Springer}
}

@article{jeong2025monotone,
  title={A monotone finite element method for an elliptic distributed optimal control problem with a convection-dominated state equation},
  author={Jeong, S. and Lee, S. and Liu, S.},
  journal={arXiv:2510.27167},
  year={2025}
}

@article{knobloch2009local,
  title={Local projection stabilization for advection-diffusion-reaction problems: {O}ne-level vs. two-level approach},
  author={Knobloch, Petr and Lube, Gert},
  journal={Applied Numerical Mathematics},
  volume={59},
  number={12},
  pages={2891--2907},
  year={2009},
  publisher={Elsevier}
}

@article{adler2023stable,
  title={A stable mimetic finite-difference method for convection-dominated diffusion equations},
  author={Adler, James H. and Cavanaugh, Casey and Hu, Xiaozhe and Huang, Andy and Trask, Nathaniel},
  journal={SIAM Journal on Scientific Computing},
  volume={45},
  number={6},
  pages={A2973--A3000},
  year={2023},
  publisher={SIAM}
}

@book{hinze2008optimization,
  title={Optimization with PDE constraints},
  author={Hinze, M. and Pinnau, R. and Ulbrich, M. and Ulbrich, S.},
  year={2008},
  publisher={Springer Science \& Business Media}
}

@article{xu1999monotone,
  title={A monotone finite element scheme for convection-diffusion equations},
  author={Xu, J. and Zikatanov, L.},
  journal={Mathematics of Computation},
  volume={68},
  number={228},
  pages={1429--1446},
  year={1999}
}

@article{ayuso2009discontinuous,
  title={Discontinuous {G}alerkin methods for advection-diffusion-reaction problems},
  author={Ayuso, B. and Marini, L. D.},
  journal={SIAM Journal on Numerical Analysis},
  volume={47},
  number={2},
  pages={1391--1420},
  year={2009},
  publisher={SIAM}
}

@article{leykekhman2012local,
  title={Local error analysis of discontinuous {G}alerkin methods for advection-dominated elliptic linear-quadratic optimal control problems},
  author={Leykekhman, D. and Heinkenschloss, M.},
  journal={SIAM Journal on Numerical Analysis},
  volume={50},
  number={4},
  pages={2012--2038},
  year={2012},
  publisher={SIAM}
}

@article{liu2024multigrid,
  title={Multigrid preconditioning for discontinuous {G}alerkin discretizations of an elliptic optimal control problem with a convection-dominated state equation},
  author={Liu, S. and Simoncini, V.},
  journal={Journal of Scientific Computing},
  volume={101},
  number={3},
  pages={79},
  year={2024},
  publisher={Springer}
}

@book{Lions,
  title={Optimal Control of Systems Governed by Partial Differential Equations},
  author={Lions, J. L.},
  year={1971},
  publisher={Springer}
}

@book{Tro,
  title={Optimal Control of Partial Differential Equations: Theory, Methods, and Applications},
  author={Tr{\"o}ltzsch, F.},
  volume={112},
  year={2010},
  publisher={American Mathematical Soc.}
}

@article{qiao2025two,
  title={Two-Scale Neural Networks for Partial Differential Equations with Small Parameters},
   author={Zhuang, Q. and Yao, C. Z. and Zhang, Z. and Karniadakis, G. E.},
  journal={Communications in Computational Physics},
  volume={38},
  number={3},
  pages={603--629},
  year={2025}
}

@article{bergounioux1992penalization,
  title={A penalization method for optimal control of elliptic problems with state constraints},
  author={Bergounioux, M.},
  journal={SIAM Journal on Control and Optimization},
  volume={30},
  number={2},
  pages={305--323},
  year={1992},
  publisher={SIAM}
}

@article{dai2025solving,
  title={Solving elliptic optimal control problems via neural networks and optimality system},
  author={Dai, Y. and Jin, B. and Sau, R. C. and Zhou, Z.},
  journal={Advances in Computational Mathematics},
  volume={51},
  number={4},
  pages={31},
  year={2025},
  publisher={Springer}
}

@book{roos2008robust,
  title={Robust numerical methods for singularly perturbed differential equations: convection-diffusion-reaction and flow problems},
  author={Roos, Hans-G{\"o}rg and Stynes, Martin and Tobiska, Lutz},
  year={2008},
  publisher={Springer}
}

@article{kingma2017adammethodstochasticoptimization,
      title={Adam: A Method for Stochastic Optimization}, 
      author={Diederik P. Kingma and Jimmy Ba},
      year={2017},
journal={arXiv:1412.6980},
      primaryClass={cs.LG},
      url={https://arxiv.org/abs/1412.6980}, 
}



\end{document}